\documentclass[11pt]{amsart}
\usepackage[utf8]{inputenc}
\usepackage[T1]{fontenc}
\usepackage{amsmath,amssymb,amsthm,mathtools}
\usepackage{booktabs}
\usepackage[margin=1.15in]{geometry}
\usepackage[dvipsnames]{xcolor}
\usepackage[colorlinks=true,linkcolor=blue!50!black,citecolor=blue!50!black,urlcolor=blue!50!black]{hyperref}

\newtheorem{theorem}{Theorem}[section]
\newtheorem{lemma}[theorem]{Lemma}
\newtheorem{proposition}[theorem]{Proposition}
\newtheorem{corollary}[theorem]{Corollary}
\newtheorem{maintheorem}{Theorem}

\theoremstyle{definition}
\newtheorem{definition}[theorem]{Definition}
\newtheorem{hypothesis}[theorem]{Hypothesis}
\newtheorem{example}[theorem]{Example}

\theoremstyle{remark}
\newtheorem{remark}[theorem]{Remark}

\newcommand{\R}{\mathbb{R}}
\newcommand{\Z}{\mathbb{Z}}
\newcommand{\T}{\mathbb{T}}
\newcommand{\Sone}{S^1}
\newcommand{\ad}{\alpha_\delta}
\newcommand{\bd}{\beta_\delta}
\newcommand{\Hbar}{\bar H}
\newcommand{\dd}{\,\mathrm{d}}
\newcommand{\avg}[1]{\langle #1\rangle}

\title[Berry phase and Hannay angle as Taylor coefficients]{Berry phase and Hannay angle\\ as Taylor coefficients of the minimal action}
\author{Karl Friedrich Siburg}
\address{Fakultät für Mathematik, TU Dortmund, Germany}
\email{karlfriedrich.siburg@tu-dortmund.de}

\begin{document}

\begin{abstract}
For a Hamiltonian system with one degree of freedom whose parameters are varied slowly around a closed loop, with slowness $\delta$, we study the minimal action $\bd$ of the resulting time-periodic system and its convex conjugate, Mather's $\alpha$-function $\ad$, both normalized per unit time. The cycle map of such a system is an area-preserving map of the cylinder which in general is not integrable. We prove that $\bd(\rho)$ and $\ad(c)$ have asymptotic expansions in $\delta$ to every finite order, locally uniformly in the rotational regime of the frozen systems, and we identify the coefficients geometrically. The zeroth coefficient of $\ad$ is the average of the frozen effective Hamiltonians, as shown by Evans and Gomes. The first coefficient is the classical Berry phase: at fixed rotation number $\rho$, the first Taylor coefficient of $\bd(\rho)$ is, modulo periods, Weinstein's Berry phase of the loop of frozen invariant curves of the class whose averaged frozen frequency is $\rho$, each curve weighted by its Mather measure. In particular, the minimal action selects the weighting in Weinstein's construction. The Hannay angle, the derivative of the Berry phase, governs the rotation numbers of all minimizing measures, including those supported on Aubry--Mather sets of the cycle map, where no invariant curve exists; and the intervals of classes on which a rotation number is locked have length $O(\delta^\infty)$. Since the expansions concern exact convex functions, no adiabatic approximation is involved. A transformation rule under loops of point transformations explains why the Berry phase is defined only modulo $2\pi c\Z$. Linear oscillators, treated by the same method, recover Berry's formula for the generalized harmonic oscillator and show that the expansions are asymptotic only.
\end{abstract}
\maketitle

\section{Introduction}

When the parameters of an integrable Hamiltonian system are moved slowly around a closed loop, the angle variables acquire, in addition to the dynamical angle $\int\omega\,\dd t$, a geometric shift, the \emph{Hannay angle} \cite{Hannay85}. Berry \cite{Berry85} showed that it is the classical counterpart of the quantum geometric phase \cite{Berry84}; the classical Berry phase and the Hannay angle are related by differentiation with respect to the action. Golin, Knauf and Marmi \cite{GKM89}, Montgomery \cite{Montgomery88} and Weinstein \cite{Weinstein90} placed these objects in a geometric framework of connections, Weinstein in particular for loops of weighted Lagrangian submanifolds.

The purpose of this paper is to show that both objects are Taylor coefficients of a single function: the minimal action $\bd$ of the slowly varied system, regarded as a function of the slowness $\delta$ at fixed rotation number, or equivalently its convex conjugate $\ad$ at fixed cohomology class. The geometric phase thus appears not as a dynamical effect that becomes visible only after an adiabatic approximation, but as the first correction coefficient of an exact convex function which is defined for every system in the class considered, integrable or not.

The idea fits into a broader picture: the Birkhoff invariants of an elliptic fixed point \cite{Siburg00} and of a positive definite invariant torus \cite{Bernard03} are Taylor coefficients of the minimal action, and for a strictly convex billiard the Marvizi--Melrose invariants \cite{MarviziMelrose} are algebraically equivalent to the Taylor coefficients of $\alpha^{2/3}$ at the boundary point $-1$, equivalently to those of the minimal action $\beta$ at $0$ \cite[\S4]{Siburg99}, \cite[Thm.~3.2.25]{Siburg04}; see \cite[Ch.~3--4]{Siburg04}. That the minimal action carries strong invariants has since been confirmed in several directions: for billiards without a strip of invariant curves at the boundary, Kaloshin and Koudjinan proved that the Taylor coefficients of a normalized $\beta$-function are invariant under $C^\infty$ conjugacy \cite{KaloshinKoudjinan}, and Huang, Kaloshin and Sorrentino showed that for generic strictly convex domains the maximal marked length spectrum, i.e.\ $\beta$ at rational rotation numbers, determines the eigendata of the Aubry--Mather periodic orbits \cite{HKS18}. Throughout we use Mather's notation: $\beta$ is the minimal action and $\alpha=\beta^*$ its convex conjugate; in \cite{Siburg99,Siburg04} these functions are denoted $\alpha$ and $\alpha^*$.

\subsection*{Setting}
Let $H_\delta(x,p,s)$ be a Hamiltonian on $T^*\T$, $\T=\R/2\pi\Z$, depending $1$-periodically on a parameter $s\in\Sone=\R/\Z$ and smoothly on $\delta\ge0$. The slowly varied system is $H_\delta(x,p,\delta t)$; one cycle of the parameter loop takes time $1/\delta$. It is a time-periodic system with one degree of freedom. Its cycle map, the time-$1/\delta$ map, is an area-preserving map of the cylinder $T^*\T$ which in general is not integrable, and Mather's theory of minimizing measures applies to it. For $c\in\R$ let $\ad(c)$ be the unique constant for which
\begin{equation}\label{eq:cell-intro}
  \delta\,\partial_s u+H_\delta(x,c+\partial_xu,s)=\ad(c)\qquad\text{on }\T\times\Sone
\end{equation}
has a Lipschitz viscosity solution. This is the effective Hamiltonian $\hat H(c,\delta)$ of Evans and Gomes \cite[(6.2)]{EvansGomes02} (see \cite{EvansGomes01} for the autonomous case); after rescaling time it is Mather's $\alpha$-function of the time-periodic system \cite{Mather91}, i.e.\ the convex conjugate of the minimal action $\bd$ (Proposition~\ref{prop:mather} and Remark~\ref{rem:mather}). Evans and Gomes proved $\lim_{\delta\to0}\ad(c)=\int_0^1\Hbar(c,s)\dd s$, where $\Hbar(\cdot,s)$ is the effective Hamiltonian of the frozen system \cite[Lemma~6.1]{EvansGomes02}, and they proposed, as a formal relation, that $\partial_\delta\ad(c)|_{\delta=0}$ equals Hannay's quantity $K(c)$ \cite[(6.19)]{EvansGomes02}.

\subsection*{Results}
For each $s$, the frozen system $H_0(\cdot,\cdot,s)$ is integrable. We call a class $c$ \emph{rotational} if every frozen system has an invariant curve $\mathcal T_{c,s}=\operatorname{graph}(c+\partial_xv(c,\cdot,s))$ of cohomology class $c$ without rest points. The rotational classes form two open half-lines $I_\pm$ (Lemma~\ref{lem:rotational} and Remark~\ref{rem:rotational}); on them the frozen curves, their invariant probability measures $\mu_{c,s}$ and their frequencies $\omega(c,s)$ depend smoothly on $(c,s)$. Let $I$ be one of these half-lines, or an open subinterval. We prove:

\begin{maintheorem}[Theorem~\ref{thm:A}]
For every $k$ there are smooth functions $\Lambda_0,\dots,\Lambda_k$ on $I$ with
\[
  \ad(c)=\sum_{j=0}^k\delta^j\Lambda_j(c)+O(\delta^{k+1})
\]
locally uniformly in $c$, where $\Lambda_0(c)=\int_0^1\Hbar(c,s)\dd s$ and
\[
  \Lambda_1(c)=K(c)=\int_0^1\!\!\int_{\T}\Bigl(\partial_s v(c,x,s)+H_1\bigl(x,c+\partial_xv(c,x,s),s\bigr)\Bigr)\dd\mu_{c,s}(x)\dd s,
\]
with $\Hbar(c,s)$ the energy of $\mathcal T_{c,s}$ and $H_1=\partial_\delta H_\delta|_{\delta=0}$.
\end{maintheorem}

By Legendre duality, Theorem~\ref{thm:A} transfers to the minimal action itself. Since $\Lambda_0'(c)=\int_0^1\omega(c,s)\dd s$ is strictly increasing (Lemma~\ref{lem:omega}), every $\rho$ in the open interval $\Lambda_0'(I)$ determines a unique class $c_0(\rho)\in I$ whose averaged frozen frequency is $\rho$.

\begin{maintheorem}[Theorem~\ref{thm:beta}]
The minimal action $\bd(\rho)$ has an asymptotic expansion in $\delta$ to every finite order, locally uniformly in $\rho\in\Lambda_0'(I)$, and, with $c_0=c_0(\rho)$,
\[
  \bd(\rho)=\beta_0(\rho)-\delta\,K(c_0)+\delta^2\Bigl(\frac{K'(c_0)^2}{2\Lambda_0''(c_0)}-\Lambda_2(c_0)\Bigr)+O(\delta^3),
\]
where $\beta_0$ is the Legendre transform of $\Lambda_0$. Every class $c$ of a minimizing measure with rotation number $\rho$ satisfies
\[
  c=c_0(\rho)-\delta\,\frac{K'(c_0)}{\Lambda_0''(c_0)}+O(\delta^2),
\]
and this expansion, too, holds to every finite order.
\end{maintheorem}

\begin{maintheorem}[Theorem~\ref{thm:B}]
For every $k$, every rotation number $\rho\in\partial\ad(c)$ of a $c$-minimizing measure satisfies, locally uniformly in $c\in I$,
\[
  \rho=\sum_{j=0}^k\delta^j\Lambda_j'(c)+O(\delta^{k+1}),\qquad\text{in particular}\qquad \rho=\int_0^1\omega(c,s)\dd s+\delta K'(c)+O(\delta^2).
\]
Moreover, on every compact subinterval of $I$, the set of classes $c$ with a prescribed rotation number $\rho_0\in\partial\ad(c)$ is an interval of length $O(\delta^\infty)$. Over one cycle, $K'(c)$ is the Hannay angle.
\end{maintheorem}

\begin{maintheorem}[Theorem~\ref{thm:C}]
If $H$ does not depend on $\delta$, then $K(c)$ equals minus Weinstein's Berry phase of the loop $s\mapsto(\mathcal T_{c,s},\mu_{c,s})$ of weighted Lagrangian submanifolds, modulo the period group $2\pi c\Z$.
\end{maintheorem}

Combining Theorems~\ref{thm:beta} and~\ref{thm:C}: \emph{at fixed rotation number $\rho$, the first Taylor coefficient of the minimal action is Weinstein's Berry phase} of the Mather-weighted loop of frozen curves of class $c_0(\rho)$, modulo $2\pi c_0(\rho)\Z$; at fixed rotation number the class is shifted by $-\delta K'(c_0)/\Lambda_0''(c_0)$, so the Hannay angle appears as this shift, rescaled by the twist $\Lambda_0''$. This is the precise sense of the title. Weinstein's construction requires a choice of weighting; Theorem~\ref{thm:C} shows that the minimal action selects the Mather measure, i.e.\ the invariant probability measure of the frozen curve. The proof uses only Weinstein's average-action formula, normalized as in his Example~7.2, and needs no sign convention for Hamiltonian vector fields beyond that normalization (Remark~\ref{rem:C}(a)). For loops in a parameter manifold, $K$ is the integral of a gauge-invariant curvature form built from the Mather densities (Remark~\ref{rem:curvature}).

Two features distinguish this approach from the adiabatic one. First, $c$ is not an approximately conserved action but the cohomology class that labels the $c$-minimizing measures of the exact time-periodic system, and the expansions concern the exact functions $\ad$ and $\bd$. The correction terms therefore arise without any adiabatic approximation; in particular, the action fluctuations which Berry and Morgan \cite{BerryMorgan96} have to take into account when they extract the geometric angle from the motion play no role here (Example~\ref{ex:rotator}). Second, $\ad$ and $\bd$ are defined for every system satisfying the standing assumptions of Section~\ref{sec:setting}, integrable or not, and Theorems~\ref{thm:B} and~\ref{thm:beta} apply to all minimizing measures, including those supported on Aubry--Mather sets of the cycle map, where no invariant curve exists.

We also prove a transformation rule for $\ad$ under loops of point transformations (Lemma~\ref{lem:invariance}): a loop with winding number $N\in\Z$ shifts $\ad(c)$ by exactly $2\pi\delta cN$. Consequently $K$ is well defined modulo $2\pi c\Z$ and the Hannay angle modulo $2\pi\Z$, which is precisely the period group in Weinstein's construction. Linear oscillators, whose invariant curves are contractible, are treated by the same corrector argument in symplectic polar coordinates (Proposition~\ref{prop:linear}); this recovers Berry's formula for the generalized harmonic oscillator. In the examples (Section~\ref{sec:examples}) the rotated rotator is solved exactly, and the rotating hoop reproduces Berry's formula; the second coefficient $\Lambda_2$ is positive in the first example and negative in the second (Remark~\ref{rem:L2}).

Here, as in the title, ``Taylor coefficients'' means the coefficients of the asymptotic expansions of $\delta\mapsto\bd(\rho)$ and $\delta\mapsto\ad(c)$ at $\delta=0$, and $O(\delta^\infty)$ means $O(\delta^N)$ for every $N$, with constants that may depend on $N$. Thus $\Lambda_1=K$ is the one-sided derivative of $\ad(c)$ at $\delta=0$, but we do not claim that $\delta\mapsto\ad(c)$ is $C^2$ there. The expansions are asymptotic only. Already for linear oscillators, $\delta\mapsto\ad(c)$ is linear on each parametric resonance interval; these intervals are invisible to the expansion (their length is $O(\delta^\infty)$), but if they accumulate at $\delta=0$ they prevent analyticity (Remark~\ref{rem:tongues}). This parallels Krikorian's result that Birkhoff normal forms generically diverge \cite{Krikorian}. For the systems of Theorem~\ref{thm:A} we do not know in general whether $\ad(c)$ is analytic in $\delta$; for the rotated rotator with analytic potential it is (Example~\ref{ex:rotator}).

We restrict to one degree of freedom, and this is not a matter of convenience. In several degrees of freedom the natural hypothesis, a smooth family of uniquely ergodic frozen invariant tori with a smoothly solvable cohomological equation, holds on no open set of classes: by convexity the frozen frequency map is a local diffeomorphism, so its image contains resonant frequencies, at which unique ergodicity fails; and even for a fixed class the frozen frequency generically crosses resonances along the loop. This is the resonance problem of adiabatic theory \cite{AKN88}.

\subsection*{Relation to previous work}
The derivations of Hannay and Berry \cite{Hannay85,Berry85} are formal: they presuppose the adiabatic invariance of the actions. In one degree of freedom the actions are adiabatic invariants for all times by KAM theory, and the existence of the Hannay angle for all initial conditions was proved by Golin \cite{Golin88}. In several degrees of freedom, Golin, Knauf and Marmi \cite[Thm.~6]{GKM89} proved, under a nondegeneracy condition and with the averaging method of Neishtadt and Bakhtin \cite{Neishtadt76,AKN88}, that the Hannay angle describes the angle shift up to an error $O(\varepsilon^b)$, for any fixed $b<1/2$, for all initial conditions outside a set of measure $O(\varepsilon^{1/2-b})$; resonances are allowed. The interpretation of the Hannay angles as the holonomy of a connection is due to Golin, Knauf and Marmi and to Montgomery \cite{GKM89,Montgomery88}, with extensions to symmetry groups \cite{MMR90} and to weighted Lagrangian submanifolds \cite{Weinstein90}; for chaotic systems the Hannay angle is replaced by ``geometric magnetism'' \cite{RobbinsBerry92}. Evans and Gomes identified $D_PK$ with the averaged correction of the angle variables under a smoothness and averaging hypothesis \cite[Thm.~5.3]{EvansGomes02} and proposed the relation $\partial_\delta\hat H(P,0)=K(P)$ formally \cite[(6.19)]{EvansGomes02}.

Accordingly, the first-order dynamical content of Theorem~\ref{thm:B} is not new for orbits on invariant curves, and in several degrees of freedom the results of \cite{GKM89} go beyond the scope of this paper. What is new here is the following.
\begin{enumerate}
\item[(i)] Berry phase and Hannay angle are Taylor coefficients of the minimal action $\bd$ at fixed rotation number and of its conjugate $\ad$ at fixed class, functions defined for every system satisfying the standing assumptions of Section~\ref{sec:setting}, integrable or not (Theorems~\ref{thm:A} and~\ref{thm:beta}). Theorem~\ref{thm:A} makes the formal relation of Evans and Gomes rigorous, to all finite orders.
\item[(ii)] Theorem~\ref{thm:B} gives the full asymptotic expansion of the rotation number of every $c$-minimizing measure of the time-periodic system, including measures supported on Aubry--Mather sets of the cycle map, and shows that rotation numbers are locked only on intervals of classes of length $O(\delta^\infty)$. It concerns these distinguished minimizing orbits, not almost all initial conditions as in \cite{Golin88,GKM89}.
\item[(iii)] The minimal action selects the Mather measure as the weighting in Weinstein's Berry phase (Theorem~\ref{thm:C}), by a proof that uses only Weinstein's average-action formula; the first Taylor coefficient of $\bd$ is this Berry phase itself.
\item[(iv)] The transformation rule (Lemma~\ref{lem:invariance}) explains the ambiguity of $K$ and of the Hannay angle by the winding of loops of coordinate changes.
\item[(v)] Linear oscillators fit into the same framework (Proposition~\ref{prop:linear}), and their resonance intervals show that the expansion is only asymptotic.
\end{enumerate}

\subsection*{Organization}
Section~\ref{sec:setting} recalls the effective Hamiltonian and proves a comparison lemma. Section~\ref{sec:A} describes the rotational regime and contains Theorem~\ref{thm:A}. Section~\ref{sec:B} contains the consequences for rotation numbers (Theorem~\ref{thm:B}) and the expansion of the minimal action (Theorem~\ref{thm:beta}). Section~\ref{sec:invariance} proves the transformation rule, Section~\ref{sec:C} the identification with Weinstein's Berry phase. Section~\ref{sec:linear} treats linear oscillators, and the final Section~\ref{sec:examples} contains two examples.

\section{The effective Hamiltonian and a comparison lemma}\label{sec:setting}

Throughout, $\T=\R/2\pi\Z$ and $\Sone=\R/\Z$. We write $f=O(\delta^\infty)$ if $f=O(\delta^N)$ for every $N\in\mathbb N$, with constants that may depend on $N$.

\begin{hypothesis}[Standing assumptions]\label{hyp:H}
For some $\delta_0>0$ the following holds.
\begin{enumerate}
\item[(S1)] The map $[0,\delta_0]\times\T\times\R\times\Sone\to\R$, $(\delta,x,p,s)\mapsto H_\delta(x,p,s)$, is $C^\infty$.
\item[(S2)] There are constants $0<\gamma\le\Gamma$ such that $\gamma\le\partial_p^2H_\delta(x,p,s)\le\Gamma$ for all $\delta,x,p,s$.
\item[(S3)] There is a constant $C$ such that $|\partial_xH_\delta(x,p,s)|\le C(1+|p|)$ for all $\delta,x,p,s$.
\end{enumerate}
\end{hypothesis}

For fixed $\delta$, (S1)--(S3) are the hypotheses (2.1)--(2.3) of \cite[\S2.1]{EvansGomes02}, periodicity in $x$ and $s$ being built into $\T\times\Sone$. Condition (S3) enters only through the a priori bound for the velocity of minimizers in \cite[proof of Thm.~2.1, Step~3]{EvansGomes02} and through the completeness of the Hamiltonian flow (see below). It excludes kinetic energies with a position-dependent metric, and it is not preserved by point transformations (Remark~\ref{rem:S3}); for this reason Lemma~\ref{lem:comparison} below is formulated under (S1) and (S2) alone. The expansion of $\ad$ (Theorem~\ref{thm:A}) and its consequences that do not refer to minimizing measures use (S3) only through the existence of a Lipschitz solution of the cell problem (Remark~\ref{rem:S3}); we therefore regard (S3) as a technical assumption. All examples of Section~\ref{sec:examples} satisfy (S1)--(S3).

We write $H_1:=\partial_\delta H_\delta|_{\delta=0}$ and $H_2:=\tfrac12\partial_\delta^2H_\delta|_{\delta=0}$. The physical system is the $1/\delta$-periodic Hamiltonian $H_\delta(x,p,\delta t)$. The case of a $\delta$-independent $H$ is the classical slowly varied system; the $\delta$-dependence is included because it arises from time-dependent coordinate changes (Section~\ref{sec:invariance}) and in examples such as the rotating hoop.

\begin{definition}\label{def:alpha}
For $\delta\in(0,\delta_0]$ and $c\in\R$, $\ad(c)$ denotes the unique real number for which the cell problem
\begin{equation}\label{eq:cell}
  \delta\,\partial_s w+H_\delta(x,c+\partial_xw,s)=\ad(c)\qquad\text{on }\T\times\Sone
\end{equation}
has a Lipschitz continuous viscosity solution $w$. Existence and uniqueness are proved in \cite[Thm.~2.1, (6.2)]{EvansGomes02}; uniqueness also follows from Lemma~\ref{lem:comparison}(b).
\end{definition}

For $u\in C^1(\T\times\Sone)$ and $c\in\R$ put
\begin{equation}\label{eq:Fu}
  F_u(x,s):=\delta\,\partial_s u(x,s)+H_\delta\bigl(x,c+\partial_xu(x,s),s\bigr).
\end{equation}

\begin{lemma}[Comparison]\label{lem:comparison}
Assume (S1) and (S2), let $\delta\in(0,\delta_0]$ and $c\in\R$, and suppose that \eqref{eq:cell}, with a constant $\lambda$ in place of $\ad(c)$, has a Lipschitz continuous viscosity solution $w$. Then:
\begin{enumerate}
\item[(a)] $\displaystyle\min_{\T\times\Sone}F_u\le\lambda\le\max_{\T\times\Sone}F_u$ for every $u\in C^1(\T\times\Sone)$;
\item[(b)] $\displaystyle\lambda=\inf_{u\in C^\infty(\T\times\Sone)}\max_{\T\times\Sone}F_u$.
\end{enumerate}
In particular $\lambda$ depends only on $H_\delta$ and $c$. If such a solution exists for every $c\in\R$, then $c\mapsto\lambda(c)$ is convex. Under Hypothesis~\ref{hyp:H}, $\lambda=\ad(c)$.
\end{lemma}

\begin{proof}
(a) Since $\T\times\Sone$ is compact, $w-u$ attains its maximum at some point $(x_0,s_0)$. There $u+\max(w-u)$ touches $w$ from above, so the subsolution property of $w$ gives $F_u(x_0,s_0)\le\lambda$, hence $\min F_u\le\lambda$. Likewise, at a minimum point of $w-u$ the supersolution property gives $F_u\ge\lambda$ there, hence $\max F_u\ge\lambda$.

(b) By (a), $\lambda\le\inf_u\max F_u$. Conversely, $w$ is differentiable almost everywhere and satisfies $\delta\,\partial_sw+H_\delta(x,c+\partial_xw,s)=\lambda$ at every point of differentiability \cite{Barles94}. Let $w_\varepsilon$ be a mollification of $w$ in $(x,s)$. The function $(a,p)\mapsto\delta a+H_\delta(x,c+p,s)$ is convex, so Jensen's inequality, together with the uniform continuity of $H_\delta$ on $\T\times\{|p|\le|c|+\operatorname{Lip}w\}\times\Sone$ (which allows us to freeze $(x,s)$ in the mollification), gives $\max F_{w_\varepsilon}\le\lambda+o(1)$ as $\varepsilon\to0$.

Convexity: let $c=tc_1+(1-t)c_2$ with $t\in[0,1]$, let $u_1,u_2\in C^\infty$ and $u:=tu_1+(1-t)u_2$. Writing $F^{c}_u$ for \eqref{eq:Fu} in the class $c$, convexity of $H_\delta$ in $p$ gives $F^{c}_u\le tF^{c_1}_{u_1}+(1-t)F^{c_2}_{u_2}$ pointwise. Taking maxima and then infima over $u_1,u_2$, (b) yields $\lambda(c)\le t\lambda(c_1)+(1-t)\lambda(c_2)$.
\end{proof}

Hypothesis~\ref{hyp:H} implies that $H_\delta$ is a Tonelli Hamiltonian. Let $L_\delta(x,v,s):=\sup_p\bigl(pv-H_\delta(x,p,s)\bigr)$ be its Lagrangian. For fixed $\delta>0$ we measure time in the slow variable $\tau=\delta t$: a curve $x(t)$ becomes $y(\tau)=x(\tau/\delta)$, and the action of the slowly varied system becomes
\[
  \int L_\delta\bigl(x,\dot x,\delta t\bigr)\dd t=\int\mathcal L_\delta\bigl(y,y',\tau\bigr)\dd\tau,\qquad
  \mathcal L_\delta(y,\eta,\tau):=\frac1\delta\,L_\delta(y,\delta\eta,\tau),
\]
a Lagrangian on $T\T\times\Sone$ which is $1$-periodic in $\tau$. Its Hamiltonian is $\tilde H:=\delta^{-1}H_\delta$. The Lagrangian $\mathcal L_\delta$ satisfies the hypotheses of Mather's theory \cite[\S1]{Mather91}, \cite[\S20]{ForniMather94}: it is $C^\infty$ and $1$-periodic in $\tau$, it is fibrewise positive definite and superlinear by (S2), and its Euler--Lagrange flow is complete. For the last point note that the Legendre transform conjugates the Euler--Lagrange flow to the Hamiltonian flow of $\tilde H$, along which $|p'|=|\partial_y\tilde H|\le\delta^{-1}C(1+|p|)$ by (S3); by Gronwall's inequality $p$, and with it the velocity $y'=\partial_p\tilde H$, stays bounded on bounded time intervals.

Let $\alpha_{\mathcal L_\delta}$ be Mather's $\alpha$-function \cite[Thm.~1]{Mather91}, \cite[Thm.~20.3]{ForniMather94},
\[
  \alpha_{\mathcal L_\delta}(c)=-\min_\mu\int\bigl(\mathcal L_\delta-c\eta\bigr)\dd\mu ,
\]
the minimum taken over all probability measures $\mu$ on $T\T\times\Sone$ with $\int\mathcal L_\delta\dd\mu<\infty$ that are invariant under the Euler--Lagrange flow of $\mathcal L_\delta$; minimizers are called \emph{$c$-minimizing}. Here $H^1(\T,\R)$ is identified with $\R$ via $c\mapsto[c\dd y]$, so that $\int\eta\dd\mu$ is Mather's rotation number of $\mu$ \cite[\S2]{Mather91}. (Mather works on the one-point compactification of $T\T\times\Sone$; invariant measures of finite action do not charge the added point.) For a $c$-minimizing measure let $\rho(\mu):=\delta\int\eta\dd\mu$ be its rotation number per unit of the physical time $t$.

\begin{proposition}\label{prop:mather}
Assume Hypothesis~\ref{hyp:H} and let $\delta\in(0,\delta_0]$.
\begin{enumerate}
\item[(a)] $\ad(c)=\delta\,\alpha_{\mathcal L_\delta}(c)$ for all $c\in\R$.
\item[(b)] $\ad$ is convex, and $\partial\ad(c)=\{\rho(\mu):\mu\text{ is $c$-minimizing}\}$.
\end{enumerate}
\end{proposition}

\begin{proof}
(a) Multiplying by $\delta$ shows that $\partial_\tau w+\tilde H(y,c+\partial_yw,\tau)=\lambda$ has a Lipschitz viscosity solution if and only if \eqref{eq:cell} does with $\ad(c)=\delta\lambda$. Hence $\lambda=\ad(c)/\delta$ is the effective Hamiltonian of $\tilde H$ in the sense of \cite[Thm.~2.1]{EvansGomes02}, and it suffices to show $\lambda=\alpha_{\mathcal L_\delta}(c)$.

\emph{$\alpha_{\mathcal L_\delta}(c)\le\lambda$.} Let $\mu$ be invariant with $\int\mathcal L_\delta\dd\mu<\infty$; then $\eta\in L^1(\mu)$ by superlinearity. Let $u\in C^1(\T\times\Sone)$. Invariance gives $\int(\partial_\tau u+\eta\,\partial_yu)\dd\mu=0$ (integrate $\frac{\dd}{\dd\tau}u$ along orbits, as in \cite[\S2]{Mather91}), and the Fenchel inequality gives $\mathcal L_\delta(y,\eta,\tau)-c\eta\ge\eta\,\partial_yu-\tilde H(y,c+\partial_yu,\tau)$. Hence
\[
  -\int\bigl(\mathcal L_\delta-c\eta\bigr)\dd\mu\le\int\bigl(\partial_\tau u+\tilde H(y,c+\partial_yu,\tau)\bigr)\dd\mu\le\frac1\delta\max F_u ,
\]
with $F_u$ as in \eqref{eq:Fu}, $s=\tau$. By Lemma~\ref{lem:comparison}(b) the infimum of the right-hand side over $u$ is $\ad(c)/\delta=\lambda$.

\emph{$\alpha_{\mathcal L_\delta}(c)\ge\lambda$.} In the proof of \cite[Thm.~2.1]{EvansGomes02} the cell problem is solved by a periodic function $g$ such that the Lax--Oleinik solution $v$ of $\partial_\tau v+\tilde H(y,c+\partial_yv,\tau)=0$, $v(\cdot,0)=g$, satisfies $v(\cdot,1)=g-\lambda$. Since the Lax--Oleinik semigroup commutes with the addition of constants and $\tilde H$ is $1$-periodic in $\tau$, $v(\cdot,k)=g-k\lambda$ for all $k\in\mathbb N$. Fix $x\in\R$ and $N\in\mathbb N$, and let $\chi\colon[0,N]\to\R$ be a minimizer in the Lax--Oleinik formula \cite[(2.9)]{EvansGomes02} for $v(x,N)$; it exists by Tonelli's theorem, and by completeness it is $C^1$ and solves the Euler--Lagrange equation \cite[\S2]{Mather91}. Then
\[
  \int_0^N\bigl(\mathcal L_\delta(\chi,\chi',\tau)-c\chi'\bigr)\dd\tau=v(x,N)-g(\chi(0))=g(x)-g(\chi(0))-N\lambda .
\]
\emph{Velocity bound.} For $k\in\{0,\dots,N-1\}$, the dynamic programming principle shows that $\chi|_{[k,k+1]}$ minimizes $\int_k^{k+1}(\mathcal L_\delta-c\xi')\dd\tau+v(\xi(k),k)$ among absolutely continuous curves $\xi$ with $\xi(k+1)=\chi(k+1)$. Since $v(\cdot,k)=g-k\lambda$ differs from $g$ by a constant and $\mathcal L_\delta$ is $1$-periodic in $\tau$, the curve $\tau\mapsto\chi(k+\tau)$, $\tau\in[0,1]$, is a minimizer in \cite[(2.9)]{EvansGomes02} for $v(\chi(k+1),1)$ with initial datum $g$. Step~3 of \cite[proof of Thm.~2.1]{EvansGomes02} bounds the velocity of every such minimizer by a constant $M$ that depends only on $\|g\|_\infty$, $c$, $\delta$ and the constants in (S2), (S3): its action is bounded in terms of $\|g\|_\infty$ \cite[(2.12)]{EvansGomes02}, so the quadratic growth of $\mathcal L_\delta$ bounds $|\chi'|$ at some time in $[k,k+1]$, and (S3) together with Gronwall's inequality bounds the momentum, hence $|\chi'|$, on all of $[k,k+1]$. Consequently $|\chi'|\le M$ on $[0,N]$, with $M$ independent of $N$ and $x$.

The empirical measures $\mu_N:=\frac1N\int_0^N\delta_{(\chi(\tau),\chi'(\tau),\tau\bmod1)}\dd\tau$ are therefore supported in the compact set $\T\times\{|\eta|\le M\}\times\Sone$, and every weak$^*$ limit $\mu$ as $N\to\infty$ is invariant under the Euler--Lagrange flow by the Kryloff--Bogoliuboff argument \cite[\S2]{Mather91}, because $\chi$ is an orbit. Since $\mathcal L_\delta$ is continuous on this compact set, dividing the last display by $N$ and passing to the limit gives $\int(\mathcal L_\delta-c\eta)\dd\mu=-\lambda$, hence $\alpha_{\mathcal L_\delta}(c)\ge\lambda$.

(b) By (a), $\ad(c)=\delta\max_\mu\bigl(c\int\eta\dd\mu-\int\mathcal L_\delta\dd\mu\bigr)$ is a maximum of affine functions of $c$, hence convex. If $\mu$ is $c'$-minimizing, then $c'$ is a subgradient of Mather's $\beta$-function $\beta=\alpha_{\mathcal L_\delta}^*$ at $\int\eta\dd\mu$ \cite[\S2, after Thm.~1]{Mather91}, \cite[p.~176]{ForniMather94}, i.e.\ $\int\eta\dd\mu\in\partial\alpha_{\mathcal L_\delta}(c')$ by Fenchel duality. Fix a compact interval $B$ containing $c$ in its interior. Since the convex function $\alpha_{\mathcal L_\delta}$ is finite, its subgradients on $B$ are bounded, so the rotation numbers of all $c'$-minimizing measures, $c'\in B$, are bounded by some $R$. The proof of \cite[Prop.~4]{Mather91} (see also \cite[Prop.~20.8]{ForniMather94}) shows that the support of a $c'$-minimizing measure is contained in $\T\times\{|\eta|\le R'\}\times\Sone$, where $R'$ depends only on such a bound $R$; it uses only that the trajectories in the support are minimizers \cite[Prop.~3]{Mather91} with asymptotic rotation number of absolute value at most $R$. Let $Q:=\T\times\{|\eta|\le R'\}\times\Sone$ and let $\mathcal Y$ be the weak$^*$ compact set of invariant probability measures supported in $Q$. For $c'\in B$ we then have $\alpha_{\mathcal L_\delta}(c')=\max_{\mu\in\mathcal Y}\varphi(c',\mu)$ with $\varphi(c',\mu):=c'\int\eta\dd\mu-\int\mathcal L_\delta\dd\mu$, which is continuous on $B\times\mathcal Y$ and affine in $c'$. Danskin's theorem \cite{Danskin67} gives that the directional derivative of $\alpha_{\mathcal L_\delta}$ at $c$ in direction $e$ is $\max\{e\int\eta\dd\mu:\mu\ c\text{-minimizing}\}$. The set of $c$-minimizing measures is convex and weak$^*$ compact, and $\mu\mapsto\int\eta\dd\mu$ is affine and continuous on $\mathcal Y$; hence $\{\int\eta\dd\mu:\mu\ c\text{-minimizing}\}$ is a compact convex set whose support function is this directional derivative. For a finite convex function the directional derivative is the support function of the subdifferential, so this set equals $\partial\alpha_{\mathcal L_\delta}(c)$. Multiplying by $\delta$ gives the claim.
\end{proof}

\begin{remark}\label{rem:mather}
Let $\beta_{\mathcal L_\delta}=\alpha_{\mathcal L_\delta}^*$ be Mather's minimal action of $\mathcal L_\delta$, i.e.\ $\beta_{\mathcal L_\delta}(h)=\min\{\int\mathcal L_\delta\dd\mu:\int\eta\dd\mu=h\}$, the minimum over invariant probability measures \cite[Thm.~1]{Mather91}. We put
\[
  \bd(\rho):=\ad^*(\rho)=\sup_{c\in\R}\bigl(\rho c-\ad(c)\bigr)=\delta\,\beta_{\mathcal L_\delta}(\rho/\delta),
\]
where the last identity follows from Proposition~\ref{prop:mather}(a). Thus $\bd(\rho)$ is the minimal action of the slowly varied system per unit physical time among invariant measures with rotation number $\rho$ per unit physical time. By Fenchel duality and Proposition~\ref{prop:mather}(b), $c\in\partial\bd(\rho)$ if and only if there is a $c$-minimizing measure with rotation number $\rho$. If a $c$-minimizing measure is supported on an invariant curve of the cycle map, $\rho(\mu)$ is the asymptotic velocity of $x$ on that curve.
\end{remark}

\section{Frozen curves and the expansion of \texorpdfstring{$\ad$}{alpha}}\label{sec:A}

\begin{hypothesis}[V]\label{hyp:V}
Let $I\subset\R$ be an open interval. There are functions $v\in C^\infty(I\times\T\times\Sone)$, $\Hbar\in C^\infty(I\times\Sone)$ and $m\in C^\infty(I\times\T\times\Sone)$ with the following properties.
\begin{enumerate}
\item[(V1)] \emph{Invariant curves:} $H_0\bigl(x,c+\partial_xv(c,x,s),s\bigr)=\Hbar(c,s)$ for all $(c,x,s)$.
\item[(V2)] \emph{Invariant density:} $m>0$, $\int_{\T}m(c,x,s)\dd x=1$, and $\partial_x(m\,X_{c,s})=0$, where
\[
  X_{c,s}(x):=\partial_pH_0\bigl(x,c+\partial_xv(c,x,s),s\bigr).
\]
\item[(V3)] \emph{Cohomological equation:} for every $f\in C^\infty(I\times\T\times\Sone)$ there is $w\in C^\infty(I\times\T\times\Sone)$ with
\[
  X_{c,s}\,\partial_xw=f-\avg f,\qquad \avg f(c,s):=\int_{\T}f(c,x,s)\,m(c,x,s)\dd x .
\]
\end{enumerate}
We write $\dd\mu_{c,s}:=m(c,\cdot,s)\dd x$.
\end{hypothesis}

By (V1), $\mathcal T_{c,s}:=\operatorname{graph}(c+\partial_xv(c,\cdot,s))$ is an invariant curve of cohomology class $c$ for the frozen Hamiltonian $H_0(\cdot,\cdot,s)$, on which the dynamics projects to the flow of $X_{c,s}$. By (V2), $\mu_{c,s}$ is invariant under this flow. By (V3) it is the only invariant probability measure: if $\nu$ is invariant, then $\int(f-\avg f)\dd\nu=\int X_{c,s}\,\partial_xw\dd\nu=0$ for all smooth $f$. Hence $\mu_{c,s}$ (lifted to $\mathcal T_{c,s}$) is the unique Mather measure of class $c$ of the frozen system \cite{CIPP98}, and $\Hbar(\cdot,s)$ is the frozen effective Hamiltonian on $I$. Moreover $X_{c,s}$ vanishes nowhere: by (V2), $mX_{c,s}$ is constant in $x$, so a zero of $X_{c,s}$ would force $X_{c,s}\equiv0$, and then (V3) fails for every $f$ that is not constant in $x$. Thus $\mathcal T_{c,s}$ is a rotational invariant curve without rest points. The following lemma shows that, conversely, such curves always give (V), and identifies the largest intervals on which (V) holds.

\begin{lemma}[The rotational regime]\label{lem:rotational}
Assume (S1) and (S2). For $s\in\Sone$ let $e_*(s):=\max_x\min_pH_0(x,p,s)$, and for $e\ge e_*(s)$ let $p_+(x,e,s)$ and $p_-(x,e,s)$ be the largest and the smallest solution $p$ of $H_0(x,p,s)=e$. Put
\[
  c_+^*:=\max_{s\in\Sone}\frac1{2\pi}\int_{\T}p_+\bigl(x,e_*(s),s\bigr)\dd x,\qquad
  c_-^*:=\min_{s\in\Sone}\frac1{2\pi}\int_{\T}p_-\bigl(x,e_*(s),s\bigr)\dd x,
\]
$I_+:=(c_+^*,\infty)$ and $I_-:=(-\infty,c_-^*)$.
\begin{enumerate}
\item[(a)] Hypothesis~\ref{hyp:V} holds on $I_+$ with
\[
  v(c,x,s)=\int_0^x\Bigl(p_+\bigl(x',\Hbar(c,s),s\bigr)-c\Bigr)\dd x',\qquad m(c,x,s)=\frac{\omega(c,s)}{2\pi X_{c,s}(x)},
\]
where $\Hbar(c,s)$ is the unique $e>e_*(s)$ with $\frac1{2\pi}\int_{\T}p_+(x,e,s)\dd x=c$ and $\omega:=\partial_c\Hbar$; moreover $X_{c,s}>0$. Likewise (V) holds on $I_-$ with $p_-$ in place of $p_+$, and $X_{c,s}<0$.
\item[(b)] Conversely, if Hypothesis~\ref{hyp:V} holds on an open interval $I$, then $I\subset I_+$ or $I\subset I_-$.
\end{enumerate}
\end{lemma}

\begin{proof}
(a) We treat $I_+$; the case of $I_-$ is analogous. By (S2), $H_0(x,\cdot,s)$ is strictly convex and grows quadratically. Hence for $e>\min_pH_0(x,p,s)$ the equation $H_0(x,p,s)=e$ has exactly two solutions $p_-<p_+$, and $\partial_pH_0(x,p_+,s)>0$. The function $p_+$ is continuous on the closed set $\{(x,e,s):e\ge e_*(s)\}$ (it is the larger root of a strictly convex function), and by the implicit function theorem it is $C^\infty$ on the open set $U:=\{(x,e,s):e>e_*(s)\}$, with $\partial_ep_+=1/\partial_pH_0(x,p_+,s)$. Consequently $c_+(e,s):=\frac1{2\pi}\int_{\T}p_+(x,e,s)\dd x$ is continuous for $e\ge e_*(s)$, $C^\infty$ for $e>e_*(s)$, and
\[
  \partial_ec_+(e,s)=\frac1{2\pi}\int_{\T}\frac{\dd x}{\partial_pH_0(x,p_+(x,e,s),s)}>0 .
\]
By the quadratic upper bound $H_0(x,p,s)\le H_0(x,0,s)+\partial_pH_0(x,0,s)\,p+\frac\Gamma2p^2$ we have $p_+\to\infty$, hence $c_+(e,s)\to\infty$, as $e\to\infty$. So for each $s$ the map $e\mapsto c_+(e,s)$ is a diffeomorphism of $(e_*(s),\infty)$ onto $(c_+(e_*(s),s),\infty)\supset I_+$, and by the implicit function theorem its inverse $\Hbar(c,s)$ is $C^\infty$ on $I_+\times\Sone$. The function $v$ is $C^\infty$, and $2\pi$-periodic in $x$ because $\int_{\T}(p_+-c)\dd x=0$; (V1) holds by construction, and $X_{c,s}(x)=\partial_pH_0(x,p_+(x,\Hbar(c,s),s),s)>0$. Differentiating $c_+(\Hbar(c,s),s)=c$ gives $\omega=\partial_c\Hbar=2\pi\bigl(\int_{\T}X_{c,s}^{-1}\dd x\bigr)^{-1}$, so $m=\omega/(2\pi X_{c,s})$ is smooth and positive, has integral $1$, and $mX_{c,s}=\omega/2\pi$ does not depend on $x$; this is (V2). For (V3) put
\[
  w(c,x,s):=\int_0^x\frac{f(c,x',s)-\avg f(c,s)}{X_{c,s}(x')}\dd x' .
\]
It is smooth, and $2\pi$-periodic in $x$ because $\int_{\T}(f-\avg f)X_{c,s}^{-1}\dd x=\frac{2\pi}\omega\int_{\T}(f-\avg f)\,m\dd x=0$.

(b) As shown before the lemma, $X_{c,s}$ vanishes nowhere, so by continuity it has a constant sign on $I\times\T\times\Sone$; say $X_{c,s}>0$. Then $c+\partial_xv(c,x,s)$ lies to the right of the minimum point of $H_0(x,\cdot,s)$, so it equals $p_+\bigl(x,\Hbar(c,s),s\bigr)$, and $\Hbar(c,s)>\min_pH_0(x,p,s)$ for every $x$, hence $\Hbar(c,s)>e_*(s)$. Integrating over $x$ gives $c=c_+(\Hbar(c,s),s)>c_+(e_*(s),s)$ for every $s$, i.e.\ $c>c_+^*$.
\end{proof}

\begin{remark}\label{rem:rotational}
Lemma~\ref{lem:rotational} also identifies the rotational classes defined pointwise in the introduction. Fix $c$, and suppose that for every $s$ the frozen system $H_0(\cdot,\cdot,s)$ has an invariant graph of class $c$ without rest points. The flow on such a graph is a single periodic orbit, so the graph lies in an energy level $e(s)$, and since the vector field has a fixed sign along it, it is the graph of $p_+(\cdot,e(s),s)$ or of $p_-(\cdot,e(s),s)$ with $e(s)>e_*(s)$. By strict monotonicity of $e\mapsto c_\pm(e,s)$, every $s$ lies in $S_+:=\{s:c>c_+(e_*(s),s)\}$ or in $S_-:=\{s:c<c_-(e_*(s),s)\}$. These sets are open, because $e_*$ and $s\mapsto c_\pm(e_*(s),s)$ are continuous, and disjoint, because $p_-\le p_+$. Since $\Sone$ is connected, $S_+=\Sone$ or $S_-=\Sone$, i.e.\ $c\in I_+$ or $c\in I_-$. Conversely, every $c\in I_\pm$ is rotational by Lemma~\ref{lem:rotational}(a).
\end{remark}

\begin{lemma}\label{lem:omega}
Under (V), $\partial_c\Hbar(c,s)=\int_{\T}X_{c,s}\dd\mu_{c,s}=:\omega(c,s)$, the frequency of the frozen curve $\mathcal T_{c,s}$, and
\[
  \partial_c\omega(c,s)=\int_{\T}\partial_p^2H_0\bigl(x,c+\partial_xv,s\bigr)\bigl(1+\partial_x\partial_cv\bigr)^2\dd\mu_{c,s}
  \ \ge\ \frac{4\pi^2\gamma}{\int_{\T}m(c,x,s)^{-1}\dd x}\ >\ 0 .
\]
\end{lemma}

\begin{proof}
Differentiating (V1) with respect to $c$ gives $\partial_c\Hbar=X_{c,s}(1+\partial_x\partial_cv)$, and differentiating once more,
\[
  \partial_c^2\Hbar=\partial_p^2H_0\,(1+\partial_x\partial_cv)^2+X_{c,s}\,\partial_x\partial_c^2v .
\]
We integrate both identities against $\mu_{c,s}$ and use $\int_{\T}X_{c,s}\,\partial_xg\,m\dd x=-\int_{\T}g\,\partial_x(mX_{c,s})\dd x=0$ for $g=\partial_cv$ and $g=\partial_c^2v$. This gives both formulas. Since $\int_{\T}(1+\partial_x\partial_cv)\dd x=2\pi$, the Cauchy--Schwarz inequality gives $4\pi^2\le\int_{\T}(1+\partial_x\partial_cv)^2m\dd x\cdot\int_{\T}m^{-1}\dd x$, and (S2) gives the lower bound.
\end{proof}

By Lemma~\ref{lem:rotational}(a), $\omega(c,s)$ is the angular frequency of the periodic frozen motion on $\mathcal T_{c,s}$: one revolution takes the time $\bigl|\int_{\T}X_{c,s}^{-1}\dd x\bigr|=2\pi/|\omega|$, and $\mu_{c,s}$ is the time average over one revolution.

\begin{definition}\label{def:K}
Under (V), for $c\in I$ we define
\begin{equation}\label{eq:K}
  K(c):=\int_0^1\!\Bigl\langle\partial_s v+H_1\bigl(\cdot,c+\partial_xv,\cdot\bigr)\Bigr\rangle(c,s)\dd s
       =\int_0^1\!\!\int_{\T}\Bigl(\partial_sv+H_1(x,c+\partial_xv,s)\Bigr)\dd\mu_{c,s}(x)\dd s.
\end{equation}
\end{definition}

The graph family determines $v$ only up to $v\mapsto v+f(c,s)$; this changes the integrand by $\partial_sf$, whose integral over the loop vanishes, so $K$ is well defined. For $\delta$-independent $H$, $K$ coincides with the quantity of Hannay and Berry in the form used by Evans and Gomes, cf.\ \cite[(5.24), (5.25), (6.19)]{EvansGomes02}.

\begin{theorem}\label{thm:A}
Assume Hypotheses~\ref{hyp:H} and~\ref{hyp:V}. For every $k\ge0$ there are $\Lambda_0,\dots,\Lambda_k\in C^\infty(I)$ such that for every compact set $\mathcal C\subset I$ there is a constant $C_k$ with
\begin{equation}\label{eq:expansion}
  \Bigl|\ad(c)-\sum_{j=0}^k\delta^j\Lambda_j(c)\Bigr|\le C_k\,\delta^{k+1}\qquad(c\in\mathcal C,\ \delta\in(0,\delta_0]).
\end{equation}
Moreover $\Lambda_0(c)=\int_0^1\Hbar(c,s)\dd s$ and $\Lambda_1=K$. The higher coefficients are given by the recursion \eqref{eq:gj}--\eqref{eq:Lj} below; for instance
\begin{equation}\label{eq:L2}
  \Lambda_2(c)=\int_0^1\Bigl\langle\partial_sw_1+\tfrac12\partial_p^2H_0\,(\partial_xw_1)^2+\partial_pH_1\,\partial_xw_1+H_2\Bigr\rangle\dd s,
\end{equation}
with all derivatives of $H$ evaluated at $(x,c+\partial_xv,s)$ and $w_1$ defined in \eqref{eq:wj}.
\end{theorem}

\begin{proof}
Fix $k$. For smooth functions $w_1,\dots,w_k$ on $I\times\T\times\Sone$ put $Y(\delta):=\sum_{i=1}^k\delta^i\partial_xw_i$ and expand the smooth function
\[
  \delta\longmapsto H_\delta\bigl(x,c+\partial_xv+Y(\delta),s\bigr)=\sum_{j=0}^k\delta^jh_j+R_k(\delta)
\]
in its Taylor polynomial at $\delta=0$ with remainder $R_k$. Then $h_0=\Hbar(c,s)$ by (V1), and for $1\le j\le k$
\[
  h_j=X_{c,s}\,\partial_xw_j+r_j ,
\]
where $r_j$ depends only on $w_1,\dots,w_{j-1}$ and on derivatives of $(\delta,p)\mapsto H_\delta$ at $(0,x,c+\partial_xv,s)$, because $w_j$ first enters at order $\delta^j$, linearly through $\partial_pH_0$. For example $r_1=H_1(x,c+\partial_xv,s)$ and $r_2=\tfrac12\partial_p^2H_0(\partial_xw_1)^2+\partial_pH_1\,\partial_xw_1+H_2$.

We now define $w_j$, $\Lambda_j$ and auxiliary functions $\tau_j$ inductively. Put $w_0:=v$, $\Lambda_0:=\int_0^1\Hbar(c,s)\dd s$, and for $j=1,\dots,k$
\begin{gather}
  g_j:=\partial_sw_{j-1}+r_j,\label{eq:gj}\\
  w_j\ \text{a smooth solution of}\ X_{c,s}\,\partial_xw_j=\avg{g_j}-g_j\quad\text{(by (V3))},\label{eq:wj}\\
  \Lambda_j(c):=\int_0^1\avg{g_j}(c,s)\dd s.\label{eq:Lj}
\end{gather}
Finally let
\[
  \tau_0(c,s):=\int_0^s\bigl(\Lambda_0(c)-\Hbar(c,\sigma)\bigr)\dd\sigma,\qquad
  \tau_j(c,s):=\int_0^s\bigl(\Lambda_j(c)-\avg{g_j}(c,\sigma)\bigr)\dd\sigma\quad(1\le j\le k).
\]
By the choice of $\Lambda_j$ these functions are $1$-periodic in $s$, and all of them are smooth in $c$. For fixed $\delta>0$ and $c\in I$ consider
\[
  u:=v+\sum_{j=1}^k\delta^jw_j+\sum_{j=0}^k\delta^{j-1}\tau_j\ \in C^\infty(\T\times\Sone).
\]
Since $\partial_x\tau_j=0$, we have $\partial_xu=\partial_xv+Y(\delta)$ and
\[
  F_u=\delta\,\partial_sv+\sum_{j=1}^k\delta^{j+1}\partial_sw_j+\sum_{j=0}^k\delta^j\partial_s\tau_j+\sum_{j=0}^k\delta^jh_j+R_k(\delta).
\]
The coefficient of $\delta^0$ is $\partial_s\tau_0+\Hbar=\Lambda_0$. For $1\le j\le k$ the coefficient of $\delta^j$ is
\[
  \partial_sw_{j-1}+\partial_s\tau_j+X_{c,s}\,\partial_xw_j+r_j
  =g_j+\bigl(\Lambda_j-\avg{g_j}\bigr)+\bigl(\avg{g_j}-g_j\bigr)=\Lambda_j .
\]
Hence
\[
  F_u-\sum_{j=0}^k\delta^j\Lambda_j=\delta^{k+1}\partial_sw_k+R_k(\delta).
\]
By Taylor's formula and compactness of $\mathcal C\times\T\times\Sone\times[0,\delta_0]$, $|R_k(\delta)|\le C'\delta^{k+1}$ uniformly; also $|\partial_sw_k|$ is bounded there. Lemma~\ref{lem:comparison} now gives \eqref{eq:expansion}. For $j=1$ we have $g_1=\partial_sv+H_1(x,c+\partial_xv,s)$, so $\Lambda_1=K$ by \eqref{eq:K}; for $j=2$ we obtain \eqref{eq:L2}.
\end{proof}

\begin{remark}\label{rem:A}
(a) The $\Lambda_j$ do not depend on the choices made (the solutions $w_j$ are unique only up to functions of $(c,s)$), since they are the coefficients of an asymptotic expansion of the intrinsically defined function $\ad$.

(b) The statement $\Lambda_0=\int\Hbar\dd s$ holds without (V) \cite[Lemma~6.1]{EvansGomes02}. Theorem~\ref{thm:A} makes the formal relation $\partial_\delta\hat H(P,0)=K(P)$ of \cite[(6.19)]{EvansGomes02} rigorous under (V), and extends it to all finite orders.

(c) In the situation of (V), the corrector $u$ is the classical two-scale expansion of averaging theory. The point of Theorem~\ref{thm:A} is that, combined with Lemma~\ref{lem:comparison}, it computes the coefficients of the exact, globally defined function $\ad$, and not of an approximate dynamics.
\end{remark}

\begin{remark}[The second coefficient]\label{rem:L2}
For $\delta$-independent $H$ we have $H_1=H_2=0$, and since $w_1$ and $m$ are $1$-periodic in $s$, an integration by parts in \eqref{eq:L2} gives
\[
  \Lambda_2(c)=\frac12\int_0^1\!\!\int_{\T}\partial_p^2H_0\,(\partial_xw_1)^2\dd\mu_{c,s}\dd s-\int_0^1\!\!\int_{\T}w_1\,\partial_sm\dd x\dd s .
\]
If the frozen densities $m(c,\cdot,s)$ do not depend on $s$, then $\Lambda_2(c)\ge0$, with equality only if $\partial_xw_1=0$, i.e.\ if $\partial_s v$ does not depend on $x$, so that the frozen curves do not move. In general, however, $\Lambda_2$ has no sign: it is positive for the rotated rotator and negative for the rotating hoop (Examples~\ref{ex:rotator} and~\ref{ex:hoop}). Note that $\Lambda_2$ is invariant under loops of point transformations (Corollary~\ref{cor:invariance}), whereas the above decomposition is not. We do not know how $\Lambda_2$ is related to the non-adiabatic corrections of Berry and Hannay \cite{BerryHannay88}.
\end{remark}

\section{Rotation numbers, the minimal action and the Hannay angle}\label{sec:B}

\begin{lemma}\label{lem:convex}
Let $f$ be convex on a neighbourhood of $[c-r,c+r]\subset\R$, let $h\in C^2$ with $|f-h|\le\varepsilon$ and $|h''|\le M$ on $[c-r,c+r]$, and assume $0<\varepsilon\le Mr^2/4$. Then $|\rho-h'(c)|\le2\sqrt{M\varepsilon}$ for every $\rho\in\partial f(c)$.
\end{lemma}

\begin{proof}
Note that $M>0$ by assumption. Let $e\in\{1,-1\}$ and $0<t\le r$. Then $f(c+te)\ge f(c)+te\rho$, hence
\[
  te\rho\le h(c+te)-h(c)+2\varepsilon\le te\,h'(c)+\tfrac M2t^2+2\varepsilon .
\]
Dividing by $t$ and choosing $t=2\sqrt{\varepsilon/M}\le r$ gives $e\bigl(\rho-h'(c)\bigr)\le2\sqrt{M\varepsilon}$. Taking $e$ equal to the sign of $\rho-h'(c)$ proves the claim.
\end{proof}

\begin{theorem}\label{thm:B}
Assume Hypotheses~\ref{hyp:H} and~\ref{hyp:V}, and let $\mathcal C\subset I$ be a compact interval.
\begin{enumerate}
\item[(a)] For every $k\ge0$ there is a constant $C_k'$ such that for all $c\in\mathcal C$, $\delta\in(0,\delta_0]$ and $\rho\in\partial\ad(c)$,
\[
  \Bigl|\rho-\sum_{j=0}^k\delta^j\Lambda_j'(c)\Bigr|\le C_k'\,\delta^{k+1}.
\]
In particular $\rho=\int_0^1\omega(c,s)\dd s+\delta K'(c)+O(\delta^2)$.
\item[(b)] Let $\kappa:=\min_{\mathcal C\times\Sone}\partial_c\omega>0$. For every $k\ge0$ there is a constant $C_k''$ such that for all $\delta\in(0,\delta_0]$, all $c_1<c_2$ in $\mathcal C$ and all $\rho_i\in\partial\ad(c_i)$,
\[
  \rho_2-\rho_1\ge\tfrac\kappa2(c_2-c_1)-C_k''\,\delta^{k+1}.
\]
In particular, for every $\rho_0\in\R$ the set $\{c\in\mathcal C:\rho_0\in\partial\ad(c)\}$ is a closed, possibly empty, interval of length at most $2C_k''\delta^{k+1}/\kappa$.
\end{enumerate}
\end{theorem}

\begin{proof}
(a) Choose $r>0$ such that $\mathcal C_r:=\{c':\operatorname{dist}(c',\mathcal C)\le r\}\subset I$. Theorem~\ref{thm:A}, with $2k+1$ in place of $k$ and $\mathcal C_r$ in place of $\mathcal C$, gives $|\ad-h_\delta|\le C\delta^{2k+2}$ on $\mathcal C_r$ for $h_\delta:=\sum_{j=0}^{2k+1}\delta^j\Lambda_j$, where we may assume $C>0$; moreover $|h_\delta''|\le M$ on $\mathcal C_r$ uniformly in $\delta\le\delta_0$, for some $M>0$. For $\delta\le\delta_1:=(Mr^2/4C)^{1/(2k+2)}$, Lemma~\ref{lem:convex} with $\varepsilon=C\delta^{2k+2}$ yields $|\rho-h_\delta'(c)|\le2\sqrt{MC}\,\delta^{k+1}$, and $h_\delta'(c)-\sum_{j=0}^k\delta^j\Lambda_j'(c)=\sum_{j=k+1}^{2k+1}\delta^j\Lambda_j'(c)$ is $O(\delta^{k+1})$ on $\mathcal C$. For $\delta\in[\delta_1,\delta_0]$ the claim holds after enlarging the constant, because $\ad$ is bounded on $\mathcal C_r$ uniformly in $\delta$ (by \eqref{eq:expansion}) and convex, so its subgradients on $\mathcal C$ are uniformly bounded. Finally $\Lambda_0'=\int_0^1\omega\dd s$ by Lemma~\ref{lem:omega}, and $\Lambda_1=K$.

(b) By Lemma~\ref{lem:omega} and compactness, $\kappa>0$, and $\Lambda_0''=\int_0^1\partial_c\omega\dd s\ge\kappa$ on $\mathcal C$. Let $L:=\max_{\mathcal C}\sum_{j=1}^k\delta_0^{j-1}|\Lambda_j''|$ (so $L=0$ if $k=0$) and $L_*:=\max\{L,\kappa/(2\delta_0)\}>0$. By (a) at $c_1$ and $c_2$ and the mean value theorem,
\[
  \rho_2-\rho_1\ge\sum_{j=0}^k\delta^j\bigl(\Lambda_j'(c_2)-\Lambda_j'(c_1)\bigr)-2C_k'\delta^{k+1}\ge(\kappa-L\delta)(c_2-c_1)-2C_k'\delta^{k+1}.
\]
For $\delta\le\kappa/(2L_*)$ we have $L\delta\le\kappa/2$, which gives the claim with $C_k''\ge2C_k'$. For $\delta>\kappa/(2L_*)$ we have $\rho_2\ge\rho_1$ by monotonicity of the subdifferential of a convex function, so the claim holds if moreover $C_k''(\kappa/2L_*)^{k+1}\ge\frac\kappa2|\mathcal C|$. Both conditions hold for $C_k'':=\max\bigl\{2C_k',\ \tfrac\kappa2|\mathcal C|\,(2L_*/\kappa)^{k+1}\bigr\}$. If $\rho_0\in\partial\ad(c_1)\cap\partial\ad(c_2)$, the inequality gives $c_2-c_1\le2C_k''\delta^{k+1}/\kappa$. Finally, by Fenchel duality, $\rho_0\in\partial\ad(c)$ if and only if $c\in\partial\bd(\rho_0)$; so the set in question is $\partial\bd(\rho_0)\cap\mathcal C$, a closed interval, possibly empty.
\end{proof}

\begin{remark}[Interpretation]\label{rem:hannay}
(a) \emph{The Hannay angle.} By Proposition~\ref{prop:mather}, every $\rho\in\partial\ad(c)$ is the average velocity, per unit time, of $x$ with respect to a $c$-minimizing measure of the slowly varied system, and conversely. The average advance of $x$ per cycle of duration $1/\delta$ is therefore
\[
  \frac1\delta\int_0^1\omega(c,s)\dd s+K'(c)+O(\delta).
\]
The first term is the dynamical angle; the second, $K'(c)$, is the Hannay angle (in radians, since $\T=\R/2\pi\Z$). No adiabatic approximation is involved: $c$ is the cohomology class of the minimizing measures of the exact system, not an approximately conserved action. The statement holds also for classes $c$ for which the cycle map has no invariant curve of class $c$. The support of a $c$-minimizing measure projects injectively onto $\T\times\Sone$, with Lipschitz inverse on its image \cite[Thm.~2]{Mather91}; in the cycle map it appears as an invariant curve, a set of periodic orbits or an Aubry--Mather set.

(b) \emph{Mode locking.} By the proof of Theorem~\ref{thm:B}(b), the set of classes with rotation number $\rho_0$ is $\partial\bd(\rho_0)\cap\mathcal C$. Thus Theorem~\ref{thm:B}(b) says that the part in $\mathcal C$ of each corner of $\bd$, i.e.\ of each subdifferential $\partial\bd(\rho_0)$, has length $O(\delta^\infty)$; applying it on a larger compact interval controls the whole corner as long as it stays inside $I$. For monotone twist maps, a corner of $\beta$ at a rational rotation number signals the absence of an invariant curve of that rotation number: $\beta$ is differentiable at a rational point if and only if there is an invariant curve consisting of periodic minimizers \cite{Mather90}. The cycle map is in general not a twist map, and we do not use an analogue of this characterization for time-periodic systems; but if it holds, the resonant rotation numbers are those with $\rho_0/(2\pi\delta)\in\mathbb Q$, i.e.\ rational rotation per cycle, and the resonance zones of the cycle map, measured in the class $c$, have width $O(\delta^\infty)$. For analytic Hamiltonians one expects exponentially small widths, in line with the exponential smallness of resonance effects in adiabatic theory \cite{AKN88}; we do not pursue this.
\end{remark}

We now turn to the minimal action itself. By Lemma~\ref{lem:omega}, $\Lambda_0''(c)=\int_0^1\partial_c\omega(c,s)\dd s>0$ on $I$, so $\Lambda_0'$ is a diffeomorphism of $I$ onto the open interval $J:=\Lambda_0'(I)$. For $\rho\in J$ let $c_0(\rho):=(\Lambda_0')^{-1}(\rho)$, the class whose averaged frozen frequency $\int_0^1\omega(c,s)\dd s$ equals $\rho$, and let
\[
  \beta_0(\rho):=\rho\,c_0(\rho)-\Lambda_0\bigl(c_0(\rho)\bigr)
\]
be the Legendre transform of $\Lambda_0$.

\begin{theorem}\label{thm:beta}
Assume Hypotheses~\ref{hyp:H} and~\ref{hyp:V}. There are functions $B_j\in C^\infty(J)$, $j\ge0$, with $B_0=\beta_0$, such that for every compact set $\mathcal J\subset J$ and every $k\ge0$ there are constants $C,\delta_*>0$ with the following properties, for all $\rho\in\mathcal J$ and $\delta\in(0,\delta_*]$.
\begin{enumerate}
\item[(a)] $\displaystyle\Bigl|\bd(\rho)-\sum_{j=0}^k\delta^jB_j(\rho)\Bigr|\le C\,\delta^{k+1}$.
\item[(b)] Every $c\in\partial\bd(\rho)$, i.e.\ every class $c$ for which some $c$-minimizing measure has rotation number $\rho$, satisfies $\displaystyle\Bigl|c-\sum_{j=0}^k\delta^jB_j'(\rho)\Bigr|\le C\,\delta^{k+1}$.
\item[(c)] With $c_0=c_0(\rho)$,
\[
  B_1(\rho)=-K(c_0),\qquad B_1'(\rho)=-\frac{K'(c_0)}{\Lambda_0''(c_0)},\qquad
  B_2(\rho)=\frac{K'(c_0)^2}{2\Lambda_0''(c_0)}-\Lambda_2(c_0).
\]
\end{enumerate}
\end{theorem}

\begin{proof}
Fix $\mathcal J$ and $k$. Choose a compact interval $\mathcal C'=[a,b]\subset I$ with $c_0(\mathcal J)\subset(a,b)$ and $\eta>0$ such that $\Lambda_0'(a)+2\eta\le\rho\le\Lambda_0'(b)-2\eta$ for all $\rho\in\mathcal J$, and let $\kappa'>0$ be a lower bound of $\Lambda_0''$ on $\mathcal C'$.

\emph{Step 1: localization.} By Theorem~\ref{thm:B}(a) with $k=0$, applied on $\mathcal C'$, there is $\delta_2>0$ such that for $\delta\le\delta_2$ every $\rho_a\in\partial\ad(a)$ satisfies $\rho_a\le\Lambda_0'(a)+\eta$ and every $\rho_b\in\partial\ad(b)$ satisfies $\rho_b\ge\Lambda_0'(b)-\eta$. Let $\rho\in\mathcal J$ and $\phi(c):=\ad(c)-\rho c$. Then $\rho_b-\rho\ge\eta$, so convexity gives $\phi(c)\ge\phi(b)+\eta(c-b)>\phi(b)$ for $c>b$, and likewise $\phi(c)>\phi(a)$ for $c<a$. Hence $\phi$ attains its minimum over $\R$, all its minimizers lie in $(a,b)$, and
\begin{equation}\label{eq:beta-local}
  \bd(\rho)=\max_{c\in\mathcal C'}\bigl(\rho c-\ad(c)\bigr),\qquad \partial\bd(\rho)=\operatorname{argmin}\phi\subset(a,b),
\end{equation}
where the second identity is Fenchel duality: $c\in\partial\bd(\rho)$ if and only if $\rho\in\partial\ad(c)$, i.e.\ if and only if $c$ minimizes $\phi$.

\emph{Step 2: comparison with a Legendre transform.} Put $A_\delta:=\sum_{j=0}^k\delta^j\Lambda_j$, defined for all real $\delta$. Since $A_0''=\Lambda_0''\ge\kappa'$ on $\mathcal C'$, there is $\delta_3>0$ such that for $|\delta|\le\delta_3$ and $\rho\in\mathcal J$ we have $A_\delta''\ge\kappa'/2$ on $\mathcal C'$, $A_\delta'(a)\le\rho-\eta$ and $A_\delta'(b)\ge\rho+\eta$. Hence $A_\delta'(c)=\rho$ has a unique solution $\hat c_\delta(\rho)\in(a,b)$, which by the implicit function theorem is smooth in $(\delta,\rho)$, with $\hat c_0=c_0$, and
\[
  \hat B_\delta(\rho):=\max_{c\in\mathcal C'}\bigl(\rho c-A_\delta(c)\bigr)=\rho\,\hat c_\delta(\rho)-A_\delta\bigl(\hat c_\delta(\rho)\bigr)
\]
is smooth in $(\delta,\rho)$ with $\partial_\rho\hat B_\delta=\hat c_\delta$, because $A_\delta'(\hat c_\delta)=\rho$. Theorem~\ref{thm:A} on $\mathcal C'$ gives $|\ad-A_\delta|\le C_k\delta^{k+1}$ there. Two functions that differ by at most $\varepsilon$ on $\mathcal C'$ have maxima over $\mathcal C'$ that differ by at most $\varepsilon$; by \eqref{eq:beta-local},
\[
  \bigl|\bd(\rho)-\hat B_\delta(\rho)\bigr|\le C_k\delta^{k+1}\qquad(\rho\in\mathcal J,\ 0<\delta\le\min\{\delta_2,\delta_3\}).
\]

\emph{Step 3: the coefficients.} Let $B_j(\rho):=\frac1{j!}\partial_\delta^j\hat B_\delta(\rho)|_{\delta=0}$. Differentiating $A_\delta'(\hat c_\delta)=\rho$ repeatedly in $\delta$ shows that $\partial_\delta^j\hat c_\delta|_{\delta=0}$, and with it $\partial_\delta^j\hat B_\delta|_{\delta=0}$, involves only $\partial_\delta^iA_\delta|_{\delta=0}=i!\,\Lambda_i$ with $i\le j$. Hence $B_j$ does not depend on $k\ge j$, nor on the choice of $\mathcal J$ and $\mathcal C'$, and so defines a smooth function on $J$. Taylor's formula in $\delta$, uniform on the compact set $\mathcal J$, gives $|\hat B_\delta-\sum_{j\le k}\delta^jB_j|\le C'\delta^{k+1}$, and (a) follows with $\delta_*:=\min\{\delta_2,\delta_3\}$. The same argument applied to $\hat c_\delta=\partial_\rho\hat B_\delta$, whose Taylor coefficients in $\delta$ are the $B_j'$, gives $|\hat c_\delta(\rho)-\sum_{j\le k}\delta^jB_j'(\rho)|\le C'\delta^{k+1}$.

For (b) let $c\in\partial\bd(\rho)$. By \eqref{eq:beta-local}, $c\in(a,b)$ and $\rho\in\partial\ad(c)$, so Theorem~\ref{thm:B}(a) on $\mathcal C'$ gives $|A_\delta'(c)-\rho|\le C_k'\delta^{k+1}$. Since $A_\delta''\ge\kappa'/2$ on $\mathcal C'$ and $A_\delta'(\hat c_\delta)=\rho$, we get $|c-\hat c_\delta(\rho)|\le2C_k'\delta^{k+1}/\kappa'$, and (b) follows.

(c) Differentiating $\hat B_\delta=\rho\hat c_\delta-A_\delta(\hat c_\delta)$ and using $A_\delta'(\hat c_\delta)=\rho$ gives $\partial_\delta\hat B_\delta=-(\partial_\delta A_\delta)(\hat c_\delta)$, hence $B_1=-\Lambda_1(c_0)=-K(c_0)$, and
\[
  \partial_\delta^2\hat B_\delta\big|_{\delta=0}=-2\Lambda_2(c_0)-K'(c_0)\,\partial_\delta\hat c_\delta\big|_{\delta=0}.
\]
Differentiating $\Lambda_0'(\hat c_\delta)+\delta K'(\hat c_\delta)+O(\delta^2)=\rho$ at $\delta=0$ gives $\partial_\delta\hat c_\delta|_{\delta=0}=-K'(c_0)/\Lambda_0''(c_0)$, which is $B_1'(\rho)$, and the formula for $B_2$ follows.
\end{proof}

\begin{remark}\label{rem:beta}
(a) \emph{Berry phase.} By Theorem~\ref{thm:C} below, for $\delta$-independent $H$ the first coefficient $B_1(\rho)=-K(c_0(\rho))$ is, modulo $2\pi c_0(\rho)\Z$, Weinstein's Berry phase of the loop $s\mapsto(\mathcal T_{c_0(\rho),s},\mu_{c_0(\rho),s})$. Thus the first Taylor coefficient of the minimal action at fixed rotation number $\rho$ is the Berry phase of the Mather-weighted frozen curves whose averaged frequency is $\rho$.

(b) \emph{Hannay angle.} Theorem~\ref{thm:beta}(b) is dual to Theorem~\ref{thm:B}(a): there the rotation number at fixed class is shifted by $\delta K'(c)$, here the class at fixed rotation number is shifted by $-\delta K'(c_0)/\Lambda_0''(c_0)$, the Hannay angle divided by the twist $\Lambda_0''=\int_0^1\partial_c\omega\dd s$. The second coefficient $B_2$ of the minimal action is not simply $-\Lambda_2$; it contains the square of the Hannay angle. Theorem~\ref{thm:beta}(b) also shows again that $\partial\bd(\rho)$ has length $O(\delta^\infty)$.

(c) \emph{Coordinate changes.} By Lemma~\ref{lem:invariance}(a) below, the conjugate of $\alpha_{H^g,\delta}$ is $\rho'\mapsto\bd(A\rho'+2\pi\delta N)$: in coordinates that move along the loop with winding number $N$, rotation numbers are shifted by $N$ full turns per cycle. Since $\beta_0'=c_0$, this shift changes $B_1$ by $2\pi Nc_0$, in accordance with Corollary~\ref{cor:invariance}.
\end{remark}

\section{Loops of coordinate changes}\label{sec:invariance}

Let $g\colon\T\times\Sone\to\T$ be smooth with $g_s:=g(\cdot,s)$ a diffeomorphism for every $s$, and write $g_s'=\partial_yg_s$. Choose a continuous lift $\hat g\colon\R\times\R\to\R$. There are $A\in\{1,-1\}$ and $N\in\Z$ such that
\[
  \hat g(y+2\pi,s)=\hat g(y,s)+2\pi A,\qquad \hat g(y,s+1)=\hat g(y,s)+2\pi N ;
\]
$A$ is the degree of $g_s$ and $N$ the \emph{winding number} of the loop $s\mapsto g_s$. (Two continuous lifts differ by a constant multiple of $2\pi$, so $N$ does not depend on the lift.) The function
\[
  h(y,s):=\hat g(y,s)-Ay-2\pi Ns
\]
is $2\pi$-periodic in $y$ and $1$-periodic in $s$.

\begin{lemma}[Transformation rule]\label{lem:invariance}
Assume Hypothesis~\ref{hyp:H}.
\begin{enumerate}
\item[(a)] \emph{Loops of point transformations.} Define
\[
  H^g_\delta(y,q,s):=H_\delta\Bigl(g_s(y),\frac q{g_s'(y)},s\Bigr)-\delta\,q\,\frac{\partial_sg_s(y)}{g_s'(y)} .
\]
Then $H^g$ satisfies (S1) and (S2) (with other constants $\gamma,\Gamma$), and for all $c\in\R$ and $\delta\in(0,\delta_0]$ the cell problem for $H^g_\delta$ in the class $Ac$ has a Lipschitz viscosity solution. Its constant, which is unique by Lemma~\ref{lem:comparison} and which we denote by $\alpha_{H^g,\delta}(Ac)$, is
\[
  \alpha_{H^g,\delta}(Ac)=\alpha_{H,\delta}(c)-2\pi\delta cN .
\]
\item[(b)] \emph{Gauge transformations.} For $f\in C^\infty(\T\times\Sone)$ let $H^f_\delta(x,q,s):=H_\delta(x,q+\partial_xf,s)+\delta\,\partial_sf$. Then $H^f$ satisfies Hypothesis~\ref{hyp:H}, and $\alpha_{H^f,\delta}=\alpha_{H,\delta}$.
\item[(c)] For $\varphi\in C^\infty(\Sone)$ one has $\alpha_{H+\varphi,\delta}=\alpha_{H,\delta}+\int_0^1\varphi\dd s$.
\end{enumerate}
\end{lemma}

The Hamiltonian $H^g$ in (a) is the one obtained from $H_\delta(x,p,\delta t)$ by the time-dependent canonical transformation $(y,q)\mapsto(g_{\delta t}(y),q/g_{\delta t}'(y))$, the cotangent lift of the loop, written in the slow variable $s=\delta t$.

\begin{proof}
(a) Property (S1) for $H^g$ is clear. Since $\partial_q^2H^g_\delta=(g_s')^{-2}\,\partial_p^2H_\delta$ and $g_s'$, $1/g_s'$ are bounded, $H^g$ satisfies (S2). Let $w$ be a Lipschitz viscosity solution of \eqref{eq:cell} for $H$ and $c$, and define
\[
  U(y,s):=w\bigl(g_s(y),s\bigr)+c\,h(y,s),
\]
a Lipschitz function on $\T\times\Sone$. For a $C^1$ function $w$ the chain rule gives, with $x=g_s(y)$ and $p=c+\partial_xw(x,s)$,
\[
  Ac+\partial_yU=g_s'(y)\,p,\qquad
  \delta\,\partial_sU=\delta\,\partial_sw(x,s)+\delta\,p\,\partial_sg_s(y)-2\pi\delta cN ,
\]
where we used $\partial_y(c\hat g)=c\,g_s'$ and $\partial_s(ch)=c\,\partial_s\hat g-2\pi cN$. Consequently
\begin{equation}\label{eq:transform}
  \delta\,\partial_sU+H^g_\delta\bigl(y,Ac+\partial_yU,s\bigr)=\delta\,\partial_sw(x,s)+H_\delta\bigl(x,c+\partial_xw,s\bigr)-2\pi\delta cN
\end{equation}
pointwise. For a viscosity solution $w$ the same holds in the viscosity sense: the map $\Psi(y,s)=(g_s(y),s)$ is a smooth diffeomorphism of $\T\times\Sone$, and $\phi$ touches $U$ from above (below) at $(y_0,s_0)$ if and only if $\psi:=\phi\circ\Psi^{-1}-c\,h\circ\Psi^{-1}$ touches $w$ from above (below) at $\Psi(y_0,s_0)$. Applying \eqref{eq:transform} to the $C^1$ function $\psi$ at that point transfers the sub- and supersolution inequalities. Hence $U$ is a Lipschitz viscosity solution of the cell problem for $H^g$ in the class $Ac$ with constant $\alpha_{H,\delta}(c)-2\pi\delta cN$. Since $A=\pm1$, every class is of the form $Ac$, and the claim follows from Lemma~\ref{lem:comparison}(b).

(b) Since $\partial_x\bigl(H_\delta(x,q+\partial_xf,s)\bigr)=\partial_xH_\delta+\partial_x^2f\,\partial_pH_\delta$ and $\partial_pH_\delta$ grows at most linearly by (S2), $H^f$ satisfies (S1)--(S3). If $w$ solves \eqref{eq:cell} for $H$, then $U:=w-f$ solves it for $H^f$ with the same constant.

(c) If $w$ solves \eqref{eq:cell} for $H$, then $w+\delta^{-1}\int_0^s(\bar\varphi-\varphi)\dd\sigma$, $\bar\varphi=\int_0^1\varphi$, solves it for $H+\varphi$ with constant $\alpha_{H,\delta}+\bar\varphi$.
\end{proof}

\begin{remark}\label{rem:S3}
Condition (S3) is in general not inherited by the Hamiltonian $H^g$ of Lemma~\ref{lem:invariance}(a): $\partial_yH^g_\delta$ contains the term $\partial_pH_\delta\cdot q\,\partial_y(1/g_s')$, which grows quadratically in $q$ unless $g_s'$ is independent of $y$, i.e.\ unless every $g_s$ is affine. This is why Lemma~\ref{lem:invariance}(a) asserts only (S1), (S2) and the solvability of the cell problem. These are all that Lemma~\ref{lem:comparison} requires, and the proof of Theorem~\ref{thm:A} uses only (S1), the existence of a Lipschitz viscosity solution of the cell problem and Lemma~\ref{lem:comparison}(a). Hence Theorem~\ref{thm:A} holds for $H^g$, which is used in Corollary~\ref{cor:invariance}. Convexity of $\alpha_{H^g,\delta}$ follows from Lemma~\ref{lem:comparison}, and with it Theorem~\ref{thm:B} and Theorem~\ref{thm:beta}(a), whose proofs use only Theorem~\ref{thm:A}, Lemma~\ref{lem:omega} and convexity. Proposition~\ref{prop:mather}, whose proof uses (S3), is not needed for $H^g$; it is needed only for the interpretation of subgradients as rotation numbers of minimizing measures.
\end{remark}

\begin{corollary}\label{cor:invariance}
Assume (V) for $H$ on $I$, and let $H^g$ be as in Lemma~\ref{lem:invariance}(a). Then $H^g$ satisfies (V) on $AI$ with
\[
  v^g(Ac,y,s)=v(c,g_s(y),s)+c\,h(y,s),\qquad m^g(Ac,y,s)=m(c,g_s(y),s)\,|g_s'(y)|,
\]
and the coefficients transform by
\[
  \Lambda^g_0(Ac)=\Lambda_0(c),\qquad K^g(Ac)=K(c)-2\pi cN,\qquad \Lambda^g_j(Ac)=\Lambda_j(c)\ (j\ge2).
\]
In particular $(K^g)'(Ac)=A\bigl(K'(c)-2\pi N\bigr)$.
\end{corollary}

\begin{proof}
The identity $\bigl(Ac+\partial_yv^g\bigr)/g_s'(y)=c+\partial_xv$ at $x=g_s(y)$ (the computation in the proof of Lemma~\ref{lem:invariance}) shows that the cotangent lift of $g_s$ maps the frozen curves of $H^g_0$ onto those of $H_0$ and conjugates the frozen dynamics; (V2) and (V3) transfer by this conjugation. The relations between the coefficients follow from Theorem~\ref{thm:A} applied to $H$ and $H^g$, Lemma~\ref{lem:invariance}(a), and uniqueness of asymptotic expansions.
\end{proof}

\begin{remark}\label{rem:periods}
Corollary~\ref{cor:invariance} shows that describing one and the same slowly varied system in coordinates that move along the loop with winding number $N$ changes, after the relabelling $c\mapsto Ac$ of classes, $K(c)$ by $2\pi cN$ and the Hannay angle, measured in the orientation of the new coordinate, by an integer multiple of $2\pi$, i.e.\ by an integer number of full turns. This is the expected ambiguity: the Hannay angle is an angle, and $2\pi c\Z$ is the group of periods of the Liouville form $p\dd x$ on $\mathcal T_{c,s}$, which is exactly the period group in Weinstein's definition of the Berry phase (Section~\ref{sec:C}).
\end{remark}

\section{The first coefficient as Weinstein's Berry phase}\label{sec:C}

We use Weinstein's conventions \cite[\S\S2, 7]{Weinstein90} for the symplectic structure: the symplectic form $\Omega=\dd x\wedge\dd p$ on $T^*\T$ (on $\R^2$ this is the form $\dd q\wedge\dd p$ of \cite[Ex.~7.2]{Weinstein90}), the Liouville form $\lambda=p\dd x$, so that $\Omega=-\dd\lambda$, and the action integral $\mathcal A(\tau)=\int_\tau\lambda$ of a loop $\tau$, so that $\mathcal A(\tau)=-\int_D\Omega$ for a spanning disc $D$ of a contractible loop. This is consistent with \cite[Ex.~7.2]{Weinstein90}, where the Berry phase of a circle revolved once counterclockwise around the origin is minus the average area swept out.

Recall Weinstein's average-action formula \cite[Thm.~7.1]{Weinstein90}: let $(L_s,\varrho_s)$, $s\in[0,1]$, be an isodrastic loop of weighted Lagrangian submanifolds obtained from $(L_0,\varrho_0)$ by a Hamiltonian isotopy $\varphi_s$ that carries $\varrho_0$ to $\varrho_s$. Then the Berry phase of the loop is the image in $G_L=\R/\Gamma_L$, where $\Gamma_L$ is the group of periods of $\lambda$ on $L_0$, of the average with respect to $\varrho_0$ of the action integrals $\mathcal A(\tau_\ell)$ of the trajectories $\tau_\ell\colon s\mapsto\varphi_s(\ell)$, $\ell\in L_0$. Weinstein derives it from his Proposition~4.4 by integration over $\ell$, and Proposition~4.4 holds with the convention $\iota_X\Omega=\dd f$ for Hamiltonian vector fields (Remark~\ref{rem:C}(a)). We take the average-action formula, which contains no Hamiltonian and is consistent with Weinstein's Example~7.2, as the operative definition of the Berry phase; then no further sign convention enters.

\begin{theorem}\label{thm:C}
Assume Hypotheses~\ref{hyp:H} and~\ref{hyp:V} with $H_\delta=H$ independent of $\delta$, and let $c\in I$. Then $s\mapsto(\mathcal T_{c,s},\mu_{c,s})$ is an isodrastic loop of positively weighted Lagrangian submanifolds, its period group is $2\pi c\Z$, and its Berry phase $B(c)$ in the sense of Weinstein satisfies
\[
  B(c)=-K(c)\mod2\pi c\Z .
\]
Consequently, for the representative $B=-K$, the Hannay angle of Theorem~\ref{thm:B} is $K'=-B'$; for any smooth choice of representative, $K'\equiv-B'$ modulo $2\pi\Z$.
\end{theorem}

\begin{proof}
The deformation $s\mapsto\mathcal T_{c,s}$ is by the exact $1$-forms $\partial_x\partial_sv\,\dd x$, hence isodrastic, and $\mu_{c,s}$ has the smooth positive density $m$ by (V2). On $\mathcal T_{c,s}$ the Liouville form restricts to $c\dd x+\dd v$, whose periods are $2\pi ck$, $k\in\Z$.

We construct symplectomorphisms carrying the weights. Let $\mathcal F(x,s):=\int_0^xm(c,x',s)\dd x'$ for $x\in\R$. Then $\partial_x\mathcal F=m>0$, $\mathcal F(x+2\pi,s)=\mathcal F(x,s)+1$, and $\mathcal F(\cdot,1)=\mathcal F(\cdot,0)$. Define $\Upsilon_s(x)$ by $\mathcal F(\Upsilon_s(x),s)=\mathcal F(x,0)$. Then $\Upsilon_s$ is a diffeomorphism of $\T$ depending smoothly on $s$, its lift satisfies $\Upsilon_0=\Upsilon_1=\mathrm{id}$, and $\Upsilon_s$ carries $\mu_{c,0}$ to $\mu_{c,s}$. Put
\[
  \varphi_s(x,p):=\Bigl(\Upsilon_s(x),\ \frac{p-c-\partial_xv(c,x,0)}{\Upsilon_s'(x)}+c+\partial_xv\bigl(c,\Upsilon_s(x),s\bigr)\Bigr).
\]
This is the composition of the fibre translation by the closed $1$-form $-(c+\partial_xv(c,\cdot,0))\dd x$, the cotangent lift of $\Upsilon_s$, and the fibre translation by $(c+\partial_xv(c,\cdot,s))\dd x$, hence a symplectomorphism of $T^*\T$. It maps $\mathcal T_{c,0}$ onto $\mathcal T_{c,s}$ and carries $\mu_{c,0}$ to $\mu_{c,s}$. Moreover $s\mapsto\varphi_s$ is a Hamiltonian isotopy with $\varphi_0=\mathrm{id}$. Indeed, at the image point $y=\Upsilon_s(x)$ the two $c$-terms combine to $c\bigl(1-1/\Upsilon_s'(x)\bigr)=c\,\partial_y\bigl(y-\Upsilon_s^{-1}(y)\bigr)$, and $y\mapsto y-\Upsilon_s^{-1}(y)$ is $2\pi$-periodic. So $\varphi_s$ is the composition of the fibre translation by the exact form $-\dd v(c,\cdot,0)$, the cotangent lift of the isotopy $\Upsilon_s$ of $\T$ starting at the identity, and the fibre translation by the exact form $\dd\bigl(c(y-\Upsilon_s^{-1}(y))+v(c,y,s)\bigr)$. Each of these is a Hamiltonian isotopy: the cotangent lift is generated by $p\,(\partial_s\Upsilon_s)\circ\Upsilon_s^{-1}$ and a translation by $\dd g_s$ by $\pm\partial_sg_s$, the sign depending on the convention.

For $\ell=(x,c+\partial_xv(c,x,0))\in\mathcal T_{c,0}$ the trajectory is $\tau_\ell(s)=\bigl(\Upsilon_s(x),c+\partial_xv(c,\Upsilon_s(x),s)\bigr)$, a closed loop because $\Upsilon_1=\Upsilon_0$ and $v$ is $1$-periodic in $s$. Since
\[
  \frac{\dd}{\dd s}\Bigl(c\,\Upsilon_s(x)+v\bigl(c,\Upsilon_s(x),s\bigr)\Bigr)=\bigl(c+\partial_xv(c,\Upsilon_s(x),s)\bigr)\,\partial_s\Upsilon_s(x)+\partial_sv\bigl(c,\Upsilon_s(x),s\bigr),
\]
and the left-hand side integrates to zero over $[0,1]$, its action integral is
\[
  \mathcal A(\tau_\ell)=\int_0^1\bigl(c+\partial_xv(c,\Upsilon_s(x),s)\bigr)\,\partial_s\Upsilon_s(x)\dd s=-\int_0^1\partial_sv\bigl(c,\Upsilon_s(x),s\bigr)\dd s .
\]
Averaging over $\mu_{c,0}$ and using that $\Upsilon_s$ carries $\mu_{c,0}$ to $\mu_{c,s}$, we obtain from \cite[Thm.~7.1]{Weinstein90}
\[
  B(c)=-\int_0^1\!\!\int_{\T}\partial_sv\bigl(c,\Upsilon_s(x),s\bigr)\dd\mu_{c,0}(x)\dd s=-\int_0^1\avg{\partial_sv}(c,s)\dd s=-K(c),
\]
since $H_1=0$.
\end{proof}

\begin{remark}\label{rem:C}
(a) \emph{Weinstein's Proposition~4.4.} The same value results from \cite[Prop.~4.4]{Weinstein90} with the time-dependent Hamiltonian $f_s(x,p):=-\partial_sv(c,x,s)$, normalized by $\tilde f_s=f_s+\avg{\partial_sv}(c,s)$, provided Hamiltonian vector fields are taken as $\dot x=\partial_pf$, $\dot p=-\partial_xf$, i.e.\ $\iota_X\Omega=\dd f$. Then $f_s$ generates the vertical flow $\dot x=0$, $\dot p=\partial_x\partial_sv$, the trajectory of $\ell$ is a closed vertical loop with $\mathcal A(\tau_\ell)=0$, and the expression of \cite[Prop.~4.4]{Weinstein90} equals $-\int_0^1\tilde f_s(\tau_\ell(s))\dd s=\int_0^1\partial_sv(c,x,s)\dd s-K(c)=-K(c)$. This convention is the one under which the expression of \cite[Prop.~4.4]{Weinstein90} is independent of $\ell$: for Weinstein's Example~7.2 a direct computation shows that it equals $-\pi(R^2+r^2)$ for every $\ell$, which agrees with Weinstein's value $-\pi R^2$ modulo the period group $\pi r^2\Z$, whereas the opposite convention gives $\pi(R^2+r^2)-2\pi|\ell|^2$. Note that the formula $H_f=-J\nabla f$ in \cite[p.~145]{Weinstein90} corresponds to the opposite convention; the proof of Theorem~\ref{thm:C} via \cite[Thm.~7.1]{Weinstein90} is not affected by this.

(b) \emph{The rotated rotator.} For the rotated rotator (Example~\ref{ex:rotator}) another family of symplectomorphisms is natural: the rotations $(x,p)\mapsto(x+2\pi s,p)$, a Hamiltonian isotopy generated by $\pm2\pi p$, carry $(\mathcal T_{c,0},\mu_{c,0})$ onto $(\mathcal T_{c,s},\mu_{c,s})$, the trajectory of $\ell=(x,p)\in\mathcal T_{c,0}$ is the loop $s\mapsto(x+2\pi s,p)$ with action $2\pi p$, and averaging over $\mu_{c,0}$ gives $2\pi\int p\dd\mu_{c,0}=2\pi\omega(c)$. Since $K(c)=2\pi(c-\omega(c))$, this representative is $B=-K+2\pi c\equiv-K$ modulo $2\pi c\Z$, and $B'=-K'+2\pi$: the two representatives differ by one full turn, as they must (Remark~\ref{rem:periods}).

(c) \emph{Hannay angle and Berry phase.} The relation $K'\equiv-B'$ is the classical relation between the Hannay angle and the Berry phase, with the cohomology class $c$ in the role of the action variable. The sign agrees with Weinstein's computation in \cite[\S9]{Weinstein90} and his Example~9.4, where, with the action variable equal to $-1/2\pi$ times the enclosed area, the Hannay curvature is $2\pi$ times the derivative of the Berry curvature with respect to the area, i.e.\ minus its derivative with respect to the action; it also agrees with Berry's formula $\Delta\theta=-\partial_I\iint W$ \cite[(15)]{Berry85}.

(d) \emph{The weighting.} Weinstein's Berry phase depends on a weighting of the Lagrangian submanifolds \cite[\S\S3, 6]{Weinstein90}. Theorem~\ref{thm:C} shows that the minimal action selects the Mather measure. In one degree of freedom $\mu_{c,s}$ is the time average over one revolution of the frozen motion, i.e.\ the Haar measure of the affine structure given by the frozen angle variable, which is the canonical weighting of toral layers in \cite[\S\S8, 9]{Weinstein90}.

(e) If $H$ depends on $\delta$, then $K$ contains the additional term $\int_0^1\avg{H_1(\cdot,c+\partial_xv,\cdot)}\dd s$, a ``magnetic'' contribution. By Lemma~\ref{lem:invariance}, a moving loop of curves in one coordinate system may become a fixed family with such a term in another; this is the situation of the rotating hoop (Example~\ref{ex:hoop}).

(f) \emph{Relation to Ahdout's fanning curves.} A curve $s\mapsto\mathcal M(s)$ of Lagrangian submanifolds is described near a point by a generating function $F(s,x)$ with $\mathcal M(s)=\operatorname{graph}\dd F_s$ \cite[\S2]{Ahdout89}; for the loop $\mathcal T_{c,s}$ one may take $F(s,x)=v(c,x,s)-v(c,x,0)$. The Hamiltonian $f_s=-\partial_sv$ that moves the curves in Remark~\ref{rem:C}(a) is thus, up to sign, the first $s$-derivative of $F$. Ahdout studies curves through a fixed point $\bar p$; under a nondegeneracy condition (``fanning'') he brings $F$ into a normal form in which $\partial_sF(0,\cdot)$ becomes a nondegenerate quadratic form \cite[Prop.~2.2]{Ahdout89}, and the symplectic invariants of the curve are the higher Taylor coefficients in $s$ \cite[Thm.~2.6]{Ahdout89}; for the geodesic flow they contain the curvature tensor \cite[Thm.~4.3]{Ahdout89}. The linearized theory, i.e.\ fanning curves in the Lagrangian Grassmannian with their horizontal curve, Jacobi endomorphism and matrix Schwarzian, is developed systematically by \'Alvarez Paiva and Dur\'an, who trace it back to Ahdout \cite[\S\S3, 6, 8.1]{AlvarezDuran09}. Theorem~\ref{thm:C} extracts a different, global invariant from the same first-order datum, namely its Mather average integrated over the loop. In this sense fanning curves provide local, pointwise invariants of a moving Lagrangian submanifold, while the minimal action records averages. We do not know whether $K$ or the higher coefficients $\Lambda_j$ can be written as averages of Ahdout's invariants; note that the loops of Theorem~\ref{thm:C} need neither pass through a fixed point nor be fanning there.
\end{remark}

\begin{remark}[Parameter spaces and curvature]\label{rem:curvature}
Suppose that the loop comes from a family on a parameter manifold $P$, i.e.\ $H(x,p,s)=G(x,p,\lambda(s))$ for a loop $\lambda\colon\Sone\to P$ and a smooth function $G$ on $T^*\T\times P$ which satisfies (S2) uniformly on compact subsets of $P$, and that $c$ is rotational for $G(\cdot,\cdot,\lambda)$ for every $\lambda\in P$. The construction in the proof of Lemma~\ref{lem:rotational}(a) is pointwise in the parameter, so it yields functions $v(c,x,\lambda)$ and $m(c,x,\lambda)$ that are smooth in $\lambda$. Consider the $1$-form
\[
  \mathcal A_c:=\int_{\T}m(c,x,\cdot)\,\dd_\lambda v(c,x,\cdot)\dd x
\]
on $P$. Since $\partial_sv=\dd_\lambda v(\lambda'(s))$, Definition~\ref{def:K} gives $K(c)=\oint_\lambda\mathcal A_c$. Replacing $v$ by $v+f(c,\lambda)$ changes $\mathcal A_c$ by $\dd_\lambda f$, because $\int_{\T}m\dd x=1$; hence the curvature
\[
  \mathcal F_c:=\dd_\lambda\mathcal A_c=\int_{\T}\dd_\lambda m(c,x,\cdot)\wedge\dd_\lambda v(c,x,\cdot)\dd x
\]
does not depend on the choice of $v$. If the loop bounds a surface $S\subset P$, then $K(c)=\iint_S\mathcal F_c$, and the Hannay angle is $K'(c)=\iint_S\partial_c\mathcal F_c$. By Theorem~\ref{thm:C}, $-\mathcal A_c$ is a connection form on the trivial bundle $P\times\R/2\pi c\Z$ whose holonomy along every loop in $P$ is Weinstein's Berry phase of the corresponding loop $\lambda\mapsto(\mathcal T_{c,\lambda},\mu_{c,\lambda})$, and $-\mathcal F_c$ is its curvature: the Berry curvature of the Mather-weighted family. For linear oscillators the analogous computation gives the surface integral of Proposition~\ref{prop:linear}(c). If the loop is not contractible in $P$, as for the rotated rotator, where $P$ is the circle of rotations, $K$ is not a surface integral, and its dependence on the loop is governed by Corollary~\ref{cor:invariance}.
\end{remark}

\section{Linear oscillators}\label{sec:linear}

The invariant curves of an oscillator are contractible and therefore invisible to Mather's theory on $T^*\T$. For linear oscillators the corrector argument nevertheless applies verbatim in symplectic polar coordinates. In these coordinates the Hamiltonian is linear in the action, so the cell problem of Definition~\ref{def:alpha} is not available; we define $\ad$ instead by the inf-max formula of Lemma~\ref{lem:comparison}(b) and prove the comparison inequalities directly.

Let $X,Y,Z\in C^\infty(\Sone)$ with $XZ-Y^2>0$, $X>0$, and consider Berry's generalized harmonic oscillator \cite{Berry85}
\[
  H(q,p,s)=\tfrac12\bigl(X(s)q^2+2Y(s)qp+Z(s)p^2\bigr),\qquad \omega(s):=\sqrt{X(s)Z(s)-Y(s)^2}.
\]
Symplectic polar coordinates $q=\sqrt{2J}\sin\theta$, $p=\sqrt{2J}\cos\theta$, $\theta\in\R/2\pi\Z$, $J>0$, satisfy $\dd q\wedge\dd p=\dd\theta\wedge\dd J$ and $p\dd q-J\dd\theta=\dd(\tfrac12J\sin2\theta)$, so actions are preserved. In these coordinates
\[
  H=J\,a(\theta,s),\qquad a(\theta,s)=X\sin^2\theta+2Y\sin\theta\cos\theta+Z\cos^2\theta>0 .
\]
The slowly varied system induces the circle flow $\dot\theta=a(\theta,\delta t)$. Let $\hat\rho_\delta$ be its rotation number per unit time, $\hat\rho_\delta=\lim_{t\to\infty}\theta(t)/t$; it exists and does not depend on the initial point (it is $\delta$ times the rotation number of the time-$1/\delta$ circle map).

\begin{definition}\label{def:alpha-lin}
For $c>0$ and $u\in C^1(\R/2\pi\Z\times\Sone)$ put $F_u:=\delta\,\partial_su+a(\theta,s)(c+\partial_\theta u)$, and
\[
  \ad(c):=\inf_{u\in C^1(\R/2\pi\Z\times\Sone)}\ \max F_u .
\]
\end{definition}

This is the characterization of Lemma~\ref{lem:comparison}(b), with the cell problem replaced by its formal analogue $\delta\,\partial_sw+a(\theta,s)(c+\partial_\theta w)=\ad(c)$ for the Hamiltonian $Ja$ in the class $c$. By Proposition~\ref{prop:linear}(a), $\ad(c)=c\,\hat\rho_\delta$. For a linear elliptic map with rotation number $\rho$ the minimal action is the indicator function of $\{\rho\}$, whose convex conjugate is $I\mapsto\rho I$; the function $c\mapsto c\,\hat\rho_\delta$ is the corresponding object for the slowly varied system, with the rotation number lifted by the flow. It is the degenerate case $b=0$ of \cite[Thm.~4.1.11]{Siburg04}.

\begin{proposition}\label{prop:linear}
\begin{enumerate}
\item[(a)] \emph{Rotation number and comparison.} For every $c>0$,
\[
  \ad(c)=c\,\hat\rho_\delta=\sup_{u\in C^1}\ \min F_u ,
\]
and consequently $\min F_u\le\ad(c)\le\max F_u$ for every $u\in C^1(\R/2\pi\Z\times\Sone)$.
\item[(b)] \emph{Expansion.} For every $k$ there are $\Theta_0,\dots,\Theta_k\in\R$ with $\hat\rho_\delta=\sum_{j=0}^k\delta^j\Theta_j+O(\delta^{k+1})$, where
\[
  \Theta_0=\int_0^1\omega(s)\dd s,\qquad
  \Theta_1=\Theta:=-\frac12\int_0^1\frac{Z(s)Y'(s)-Y(s)Z'(s)}{Z(s)\,\omega(s)}\dd s .
\]
\item[(c)] Consequently $\ad(c)=c\bigl(\Theta_0+\delta\Theta+O(\delta^2)\bigr)$, so $K(c)=\Theta c$, and the Hannay angle $\partial_cK=\Theta$ is independent of the action. By Stokes' theorem, for any surface $S$ in $\{XZ>Y^2\}$ bounded by the loop,
\[
  \Theta=\frac14\iint_S\frac{X\dd Y\wedge\dd Z+Y\dd Z\wedge\dd X+Z\dd X\wedge\dd Y}{(XZ-Y^2)^{3/2}},
\]
which is Berry's result $\Delta\theta=-\partial_I\iint_SW$ with his angle $2$-form $W$ \cite[(15), (38)]{Berry85}; equivalently, our $1$-form is $-\tfrac12(Z/\omega)\dd(Y/Z)$, cf.\ \cite[(37)]{Berry85}.
\end{enumerate}
\end{proposition}

\begin{proof}
(a) Let $\mathcal X:=a\,\partial_\theta+\delta\,\partial_s$ on $\R/2\pi\Z\times\Sone$ and $\Phi_t$ its flow; the $\theta$-component of $\Phi_t(\theta_0,s_0)$ solves $\dot\theta=a(\theta,s_0+\delta t)$. Then $F_u=c\,a+\mathcal Xu$, so along an orbit $F_u(\Phi_t)=c\,a(\Phi_t)+\frac{\dd}{\dd t}u(\Phi_t)$. Averaging over $[0,T]$ and letting $T\to\infty$ gives $c\hat\rho_\delta$, since $u$ is bounded, and the average lies between $\min F_u$ and $\max F_u$. Hence $\sup_u\min F_u\le c\hat\rho_\delta\le\inf_u\max F_u$.

For the converse inequalities let $\tilde\theta(t)$ denote a lift to $\R$ of the $\theta$-component of $\Phi_t(\theta_0,s_0)$. We claim that $|\tilde\theta(T)-\tilde\theta(0)-T\hat\rho_\delta|\le C_\delta$ for all $T\ge0$ and all initial points. Indeed, for $s_0=0$ and $T=k/\delta$ this is the standard estimate $|F^k(\vartheta)-\vartheta-k\rho(F)|<2\pi$ for a lift $F$ of the time-$1/\delta$ map, a circle diffeomorphism with rotation number $\rho(F)=\hat\rho_\delta/\delta$. A general orbit reaches $s\in\Z$ after a time at most $1/\delta$, and between two such times $\tilde\theta$ changes by at most $\max a/\delta$; this gives the claim with $C_\delta:=2\pi+2(\max a+\hat\rho_\delta)/\delta$. Now put
\[
  u_T:=\frac cT\int_0^T\!\!\int_0^t a\circ\Phi_r\dd r\dd t\in C^\infty(\R/2\pi\Z\times\Sone).
\]
Since $\mathcal X(a\circ\Phi_r)=\frac{\dd}{\dd r}(a\circ\Phi_r)$, we have $\mathcal Xu_T=\frac cT\int_0^T(a\circ\Phi_t-a)\dd t$, hence
\[
  F_{u_T}=\frac cT\int_0^Ta\circ\Phi_t\dd t=c\,\frac{\tilde\theta(T)-\tilde\theta(0)}T ,
\]
evaluated along the orbit of the respective point. By the claim, $|F_{u_T}-c\hat\rho_\delta|\le cC_\delta/T$ everywhere, so $\inf_u\max F_u\le c\hat\rho_\delta\le\sup_u\min F_u$.

(b) \emph{Frozen data.} Let $\varphi(\theta,s)$ be the frozen angle variable, i.e.\ $\dot\varphi=\omega(s)$ along the frozen flow, normalized to be a degree-one circle diffeomorphism in $\theta$; then $\partial_\theta\varphi=\omega/a$. The frozen curve of action $c$ is $J=c\,\omega/a=c(1+\partial_\theta\Xi)$ with $\Xi:=\varphi-\theta$, and its invariant probability measure is $\dd\mu_s=\dd\varphi/2\pi=(\omega/2\pi a)\dd\theta$. For smooth $f$ the equation $a\,\partial_\theta w=\avg f-f$, $\avg f=\int f\dd\mu_s$, has the smooth periodic solution $w=\int_0^\theta(\avg f-f)/a\dd\theta'$, exactly as in the proof of Lemma~\ref{lem:rotational}(a).

\emph{Corrector.} Put $w_0:=\Xi$, define recursively $g_j:=\partial_sw_{j-1}$, $w_j$ by $a\,\partial_\theta w_j=\avg{g_j}-g_j$, $\Theta_j:=\int_0^1\avg{g_j}\dd s$ ($j\ge1$), $\Theta_0:=\int\omega\dd s$, and $\tau_j$ as in the proof of Theorem~\ref{thm:A} (with $\Hbar$ replaced by $\omega$). For $u:=c\bigl(\Xi+\sum_{j=1}^k\delta^jw_j+\sum_{j=0}^k\delta^{j-1}\tau_j\bigr)$ one finds, because $F_u$ is affine in $u$ and no Taylor remainder occurs,
\[
  F_u=c\sum_{j=0}^k\delta^j\Theta_j+\delta^{k+1}c\,\partial_sw_k
\]
exactly. Part (a) gives (b), with $\Theta_1=\int_0^1\avg{\partial_s\Xi}\dd s$.

\emph{Computation of $\Theta_1$.} The frozen action-angle variables are
\[
  q=\sqrt{2IZ/\omega}\,\sin\varphi,\qquad p=\sqrt{2I/(Z\omega)}\,(\omega\cos\varphi-Y\sin\varphi),
\]
as one checks directly ($H=\omega I$, $\dd q\wedge\dd p=\dd\varphi\wedge\dd I$). Hence $\varphi(\theta,s)$ is defined by $G(\varphi,\theta,s):=Z\sin\varphi\cos\theta-(\omega\cos\varphi-Y\sin\varphi)\sin\theta=0$. Since $\Xi=\varphi-\theta$, we have $\partial_s\Xi=\partial_s\varphi|_\theta=-\partial_sG/\partial_\varphi G$. On the zero set, $(\sin\theta,\cos\theta)$ is proportional to $(Z\sin\varphi,\ \omega\cos\varphi-Y\sin\varphi)$; inserting this,
\[
  \partial_\varphi G\propto Z\omega,\qquad
  \partial_sG\propto Z'\sin\varphi(\omega\cos\varphi-Y\sin\varphi)-Z\sin\varphi(\omega'\cos\varphi-Y'\sin\varphi),
\]
with the same factor. Averaging over $\varphi$ with respect to $\dd\varphi/2\pi$ gives $\avg{\partial_s\Xi}=-(ZY'-YZ')/(2Z\omega)$, hence the formula for $\Theta$.

(c) The first statements follow from (b). For the surface integral one checks that $\dd$ of the $1$-form $-(Z\dd Y-Y\dd Z)/(2Z\omega)$ on $\{XZ>Y^2\}$ equals the integrand.
\end{proof}

\begin{remark}\label{rem:linear}
(a) Hannay's form of the angle, $-\partial_I\avg{p\,\partial_sq}$ in frozen action-angle variables, differs from the integrand $\avg{\partial_s\Xi}$ by the exact derivative $\frac{\dd}{\dd s}\bigl(Y/2\omega\bigr)$, so both give the same $\Theta$ over a loop.

(b) In polar coordinates the frozen action-angle variables $(\varphi,I)$ are the cotangent lift of the circle diffeomorphism $\theta\mapsto\varphi(\theta,s)$, since $I=J/\partial_\theta\varphi$. This is special to linear systems; for nonlinear oscillators the change to frozen action-angle variables is not a point transformation.
\end{remark}

\begin{remark}[Resonance intervals]\label{rem:tongues}
If the monodromy of the linear $1/\delta$-periodic system is hyperbolic (parametric resonance), the rotation per cycle is an integer multiple of $\pi$, so $\hat\rho_\delta=\nu\pi\delta$ is linear on each such $\delta$-interval, with $\nu\in\Z$. These intervals have length $O(\delta^\infty)$. Indeed, let $[\delta_a,\delta_b]$ be such an interval, fix $k$, and put $P_k(\delta):=\sum_{j=0}^k\delta^j\Theta_j$. By (b), $|\nu\pi\delta-P_k(\delta)|\le C_k\delta^{k+1}$ at $\delta=\delta_a$ and at $\delta=\delta_b$. The same estimate at $\delta_a$ gives $\nu\pi=(\Theta_0+O(\delta_a))/\delta_a\ge\Theta_0/(2\delta_b)$ for small $\delta_b$, where $\Theta_0>0$, whereas $|P_k'|\le M_k$ on $[0,\delta_0]$. Subtracting the two estimates and applying the mean value theorem to $\nu\pi\delta-P_k(\delta)$ yields $(\nu\pi-M_k)(\delta_b-\delta_a)\le2C_k\delta_b^{k+1}$, hence $\delta_b-\delta_a\le(8C_k/\Theta_0)\,\delta_b^{k+2}$ once $M_k\le\Theta_0/(4\delta_b)$. The intervals are therefore invisible to the expansion; but since $\nu$ changes from one interval to the next, $\delta\mapsto\ad(c)$ is not real-analytic whenever such intervals accumulate at $\delta=0$, and generically $\hat\rho_\delta$ has square-root singularities at their endpoints, where the trace of the monodromy passes through $\pm2$. In that case the expansion of Proposition~\ref{prop:linear} is asymptotic only, in line with \cite{Krikorian}. We do not know whether an analogous phenomenon occurs in the setting of Theorem~\ref{thm:A}.
\end{remark}

\section{Examples}\label{sec:examples}

\begin{example}[Rotated rotator]\label{ex:rotator}
Let $V$ be smooth and $2\pi$-periodic and $H(\xi,p,s)=\tfrac12p^2+V(\xi-2\pi s)$ on $T^*(\R/2\pi\Z)$: a potential rotated once per cycle \cite{Berry85,BerryMorgan96}. The same example appears in \cite[\S2, Remark~1]{GKM89} in the low-energy regime, where the frozen curves are contractible; here we consider the rotational regime. Let $\Hbar$ be the effective Hamiltonian of $\tfrac12p^2+V$ and $\omega=\Hbar'$. Here $e_*=\max V$ and $p_\pm=\pm\sqrt{2(e-V)}$, so the rotational regime of Lemma~\ref{lem:rotational} is $|c|>c_0:=\frac1{2\pi}\int_0^{2\pi}\sqrt{2(\max V-V)}\dd\xi$, and (V) holds on $I=(c_0,\infty)$.

Apply Lemma~\ref{lem:invariance}(a) with $g_s(y)=y+2\pi s$, so $A=1$, $N=1$, $h=0$. Then
\[
  H^g_\delta(y,q)=\tfrac12q^2+V(y)-2\pi\delta q=\tfrac12(q-2\pi\delta)^2+V(y)-2\pi^2\delta^2
\]
is autonomous, so its cell problem has $s$-independent solutions and $\alpha_{H^g,\delta}(c)=\Hbar(c-2\pi\delta)-2\pi^2\delta^2$. Hence, \emph{exactly},
\[
  \ad(c)=\Hbar(c-2\pi\delta)-2\pi^2\delta^2+2\pi\delta c,
\]
and
\[
  K(c)=2\pi\bigl(c-\omega(c)\bigr),\qquad \partial_cK(c)=2\pi\bigl(1-\omega'(c)\bigr),\qquad \Lambda_2(c)=2\pi^2\bigl(\omega'(c)-1\bigr).
\]
The Hannay angle $\partial_cK$ is the one found by Berry and Morgan \cite[(12)]{BerryMorgan96}, who obtain it from the average over the frozen curve. They also show that the action fluctuations of order $\delta$, which accumulate over the time $1/\delta$, must be taken into account when the geometric angle is extracted from the motion. Here the value follows from the exact formula for $\ad$, without any approximation. The same value of $K$ results from Definition~\ref{def:K} directly: with $v(c,\xi,s)=v_0(c,\xi-2\pi s)$ and $P:=c+v_0'$ one has $\dd\mu=\frac{\omega}{2\pi}\frac{\dd\xi}{P}$ and $\avg{v_0'}=\omega-c$, so $\avg{\partial_sv}=2\pi(c-\omega)$.

Since $c$ is the arithmetic mean and $\omega(c)$ the harmonic mean of $P$ over the circle, $K(c)>0$ for $c>c_0$ unless $V$ is constant, so $\delta\mapsto\ad(c)$ increases for small $\delta$. Moreover, writing $\overline{f}$ for the mean of $f$ over the circle, $c=\overline P$ with $P=\sqrt{2(e-V)}$ and $e=\Hbar(c)$ gives $\dd c/\dd e=\overline{P^{-1}}$ and $\dd^2c/\dd e^2=-\overline{P^{-3}}$, hence
\[
  \omega'(c)=\frac{\overline{P^{-3}}}{\bigl(\overline{P^{-1}}\bigr)^3}\ \ge\ 1
\]
by the inequality between power means, with equality only for constant $V$. Thus $\Lambda_2>0$ on $I$ unless $V$ is constant (Remark~\ref{rem:L2}). Since $H$ does not depend on $\delta$, the exact formula holds for all $\delta>0$. Let $\bar V$ and $\operatorname{Var}V$ be mean and variance of $V$ over the circle. Expanding $c=\frac1{2\pi}\int_0^{2\pi}\sqrt{2(e-V)}\dd\xi$ for large $e$ gives $\Hbar(p)=\frac12p^2+\bar V+\operatorname{Var}V/(2p^2)+O(p^{-4})$ as $|p|\to\infty$, hence $\ad(c)\to\frac12c^2+\bar V$ as $\delta\to\infty$. By Jensen's inequality, $c\le\sqrt{2(\Hbar(c)-\bar V)}$, i.e.\ $\Hbar(c)\ge\frac12c^2+\bar V$, with equality only for constant $V$. Thus $\ad(c)$ exceeds $\Hbar(c)=\lim_{\delta\to0}\ad(c)$ for small $\delta$ and falls below it for large $\delta$; it is not monotone in $\delta$ unless $V$ is constant. Compare the monotonicity in the time-frequency parameter found by Liu \cite{Liu25} in a setting without cohomology class.
\end{example}

\begin{example}[Rotating hoop]\label{ex:hoop}
A bead slides without friction on a closed plane curve $\gamma$ of length $\ell$ enclosing the area $\mathcal S$, parametrized counterclockwise by arc length $\sigma$; the curve is rotated once per cycle about the origin with angular velocity $\varpi=2\pi\delta$ \cite{Hannay85,Berry85}. In the inertial frame the Lagrangian is $\tfrac12|\gamma'\dot\sigma+\varpi R\gamma|^2=\tfrac12\dot\sigma^2+\varpi b(\sigma)\dot\sigma+\tfrac12\varpi^2|\gamma|^2$ with $R$ the rotation by $\pi/2$ and $b=\gamma'\cdot R\gamma$. In the variable $x=k\sigma\in\R/2\pi\Z$, $k=2\pi/\ell$,
\[
  H_\delta(x,p)=\frac{k^2}2\Bigl(p-\frac{2\pi\delta}k\,b\Bigr)^2-2\pi^2\delta^2|\gamma|^2 .
\]
This is autonomous and depends on $\delta$; the frozen system is $H_0=\tfrac12k^2p^2$ with invariant curves $p=c$, $v=0$ and uniform $\mu$. Since $H_1=-2\pi k\,b\,p$ and $\bar b:=\frac1{2\pi}\int b\dd x=\frac1\ell\oint(\gamma_1\dd\gamma_2-\gamma_2\dd\gamma_1)=2\mathcal S/\ell$,
\[
  K(c)=-\frac{8\pi^2\mathcal S}{\ell^2}\,c,\qquad \partial_cK=-\frac{8\pi^2\mathcal S}{\ell^2}.
\]
This is Berry's formula \cite[(55)]{Berry85}; for a circle ($\mathcal S=\ell^2/4\pi$) it gives $-2\pi$, the slippage of the reference point \cite[(56)]{Berry85}. For the second coefficient, \eqref{eq:wj} gives $\partial_xw_1=2\pi(b-\bar b)/k$, and $H_2=2\pi^2(b^2-|\gamma|^2)$; inserting this into \eqref{eq:L2} yields
\[
  \Lambda_2(c)=2\pi^2\Bigl(\frac{4\mathcal S^2}{\ell^2}-\overline{|\gamma|^2}\Bigr),
\]
where $\overline{|\gamma|^2}$ is the mean of $|\gamma|^2$ with respect to arc length. The same value follows from an expansion of the exact effective Hamiltonian of the autonomous system $H_\delta$. Since $|b|\le|\gamma|$, we have $\bar b^2\le\overline{b^2}\le\overline{|\gamma|^2}$, so $\Lambda_2\le0$, with equality if and only if $\gamma$ is a circle centred at the origin.
\end{example}

\section*{Declaration of generative AI and AI-assisted technologies in the manuscript preparation process}

During the preparation of this work the author used Claude (Anthropic) in order to check proofs and computations, and to assist with editing the text. After using this tool, the author reviewed and edited the content as needed, verified all mathematical statements and proofs independently, and takes full responsibility for the content of the publication.


\end{document}